\documentclass[a4paper,12pt]{article}
\usepackage{amsmath}
\usepackage[mathscr]{eucal}
\usepackage{amssymb}
\usepackage{latexsym}
\usepackage{amsthm}
\usepackage{amscd}
\theoremstyle{plain}
\newtheorem{theorem}{Theorem}
\newtheorem{corollary}{Corollary}[section]
\newtheorem{lemma}{Lemma}[section]
\newtheorem{proposition}{Proposition}

\theoremstyle{remark}
\newtheorem{remark}{Remark}
\theoremstyle{definition}

\usepackage[dvips]{graphicx}
\usepackage{psfrag}
\begin{document}
\title{ Surfaces with $c_1^2 =9$ and $\chi =5$  whose \\ 
canonical  classes are divisible by $3$\footnotetext{The author acknowledges the support by JSPS Grant-in-Aid Scientific Research(C) 15K04825.}}
\author{Masaaki Murakami}
\date{}
\maketitle
\begin{abstract}
We shall study minimal complex   
surfaces with $c^2 = 9$ and $\chi=5$ whose canonical 
classes are divisible by $3$ in the integral cohomology groups,  
where $c_1^2$ and $\chi$ denote the first Chern number of an 
algebraic surface 
and the Euler characteristic of the structure sheaf, 
respectively. 
The main results are a structure theorem for such surfaces, 
the unirationality of the moduli space, and a description of 
the behavior of the canonical map.  
As a byproduct, we shall also rule out a certain case 
mentioned in a paper by Ciliberto--Francia--Mendes Lopes. 
Since the irregularity $q$ vanishes for our surfaces, 
our surfaces have geometric genus $p_g = 4$.  
\end{abstract}

\section{Introduction}  

When one wants to study the behavior of canonical maps of 
algebraic surfaces, 
surfaces of general type with 
$p_g = 4$ are in a sense the most primitive objects,  
since their canonical images are  in 
most cases hypersurfaces  
of the $3$-dimensional projective space $\mathbb{P}^3$. 
Partly for such reasons, these surfaces have attracted
many algebraic geometers, 
even from the time of classical Italian school.

After Noether and Enriques studied the case $c_1^2 =4$, 
surfaces with $p_g=4$ have been studied from various 
view points 
(e.g, \cite{catbabbage}, \cite{cilcansfpg4}, \cite{kodordsings}).    
As for the classification, 
Horikawa and Bauer completed that   
for the surfaces of cases $4 \leq c_1^2 \leq 7$  
(\cite{smallc1-1}, \cite{quintic}, \cite{smallc1-3}, \cite{pg4c7bauer}). 
Complete classification of the surfaces of case $c_1^2 = 8$
seems not completely out of reach, but for the moment, only partial 
classifications and several examples are known 
(e.g.,\cite{caninvpg4c8}, \cite{cilcansfpg4}, \cite{onbicanonicalmaps}). 
We also notice that even though the surfaces 
have been classified for the case $c_1^2 =6$, 
the number of the irreducible components of the 
moduli space remains unknown even after \cite{bcpmodulic6}.

Among such works, the results most connected to 
the present paper are those on even surfaces for the case $c_1^2 = 8$. 
Recall that an algebraic surface is said to be even 
if its canonical class is divisible by $2$. 
In \cite{oliverioevenpg4c8}, Oliverio studied regular even surfaces 
of case $c_1^2 =8$, and showed that if $S$ is a surface of 
this class with base point free canonical system, then 
its canonical model is a $(6, 6)$-complete intersection 
in the weighted projective space $\mathbb{P} (1,1,2,3,3)$. 
He also showed that these surfaces fill up an open dense 
subset of a $35$-dimensional irreducible 
component $\mathcal{M}_{\mathcal{F}}$ of the moduli space 
$\mathcal{M}_{8, 4, 0}^{\mathrm{ev}}$ of even regular surfaces 
of case $c_1^2=8$. 
Though this paper \cite{oliverioevenpg4c8} studied these 
surfaces only under the condition 
that the canonical systems are base point free, 
Catanese, Liu, and Pignatelli later in \cite{catlipign'} classified  
all even regular surfaces  with $c_1^2=8$ and $p_g = 4$ and 
showed that $\mathcal{M}_{8, 4, 0}^{\mathrm{ev}}$ consists exactly 
of two irreducible components $\mathcal{M}_{\mathcal{F}}$ and 
$\mathcal{M}_{\mathcal{E}}$, both of dimension $35$ and intersecting 
each other in codimension one locus.   

In this paper, we go one step up, and study regular surfaces 
of case $c_1^2 = 9$ with canonical classes divisible by $3$. 
We shall prove three theorems. 
Our first theorem asserts that any surface of this class has the 
canonical model isomorphic to  a  $(6, 10)$-complete intersection 
of the weighted projective space $\mathbb{P} (1, 2, 2, 3, 5)$ 
(Theorem \ref{thm:structurethm}). 
Our second theorem asserts that the moduli space of our 
surfaces is unirational of dimension $34$, hence also the uniqueness
of the diffeomorphic tpye of our surfaces (Theorem \ref{thm:modulisp}). 
Our third theorem asserts that the canonical map 
$\varPhi_{|K|}$ of a surface of this class is either 
birational onto a singular sextic or generically 
two-to-one onto a cubic surface (Theorem \ref{thm:canonicalimage}).    
The surfaces with birational $\varPhi_{|K|}$ 
and those with generically two-to-one $\varPhi_{|K|}$
form an open dense subset
and a $33$-dimensional locus, respectively,  
in $\mathcal{M}$.    
  
Possibility of the existence of surfaces 
with $c_1^2 = 9$ and $p_g= 4$ and with 
canonical classes divisible by $3$ has already been mentioned 
in \cite[(ii), Proposition 1.7]{onbicanonicalmaps}, 
though for the case of canonical map composite 
with a pencil.    
In fact, the construction of examples of Case (ii) above 
was one of the motivations for our work. 
In the course of the proof of our Theorem \ref{thm:structurethm}, 
however, we shall 
show that this Case (ii) never occurs, even for the case of 
positive irregularity (Proposition \ref{prop:nevoccurcfml}).  
This sharpens their Proposition 1.7 slightly.

Let $L$ be a divisor linearly equivalent to the canonical 
divisor of our surface. 
Our strategy of the first part is 
to study the map $\varPhi_{|2L|}$ to compute the dimensions 
of some cohomology groups, where $\varPhi_{|2L|}$ is the map 
associated to the linear system $|2L|$.
Although the main tools for this part are classical ones, e.g., 
the double cover technique,   
a result by the author given in \cite{remarks'''''}
on the torsion groups of 
surfaces with $c_1^2 = 2 \chi - 1$ 
is also used to rule out some cases. 
Then we divide our argument into two cases depending on 
whether $\varPhi_{|2L|}$ is composite with a pencil or not,   
and study each case.  
For the case where $\varPhi_{|2L|}$ is composite with a pencil,  
it turns out that we are in Case (ii) of 
\cite[Proposition 1.7]{onbicanonicalmaps}. 
Using results in \cite{onbicanonicalmaps} and applying to 
$\varPhi_{|2L|}$ the structure theorem for genus $3$ fibrations 
given in \cite{fibrationsI'}, 
we shall rule out this case. 
For the case where $\varPhi_{|2L|}$ is not composite with a 
pencil, we shall study the semicanonical ring 
$R = \bigoplus_{n=0}^{\infty} H^0 (\mathcal{O} (nL))$. 
Using arguments similar to those in \cite{globalpg1k2=1}, we shall 
find generators of the ring $R$ and relations among them, 
which gives us the structure theorem. 
As for the results on the moduli space and 
the canonical maps,  
we shall prove them using this structure theorem.   
In addition to the theorems stated above, 
we shall also give a double cover 
description of our surfaces with $\deg \varPhi_{|K|} =2$ 
(Proposition \ref{prop:branchdivofphiK}).

After all the main results of the present paper were obtained, 
Kazuhiro Konno pointed out the normality of the canonical images of our 
surfaces of case $\deg \varPhi_{|K|} = 1$ 
(Proposition \ref{prop:normalityim}).  
As informed to the author by him, our surfaces therefor give 
one of the missing examples of the list given in Konno's work 
\cite{normalcankonno} on normal canonical surfaces with $p_g = 4$.   

{\sc Acknowledgment}  
 
The author expresses his gratitude to Prof.\,Kazuhiro Konno,    
who kindly gave him the comment on the normality of the canonical 
image (Proposition \ref{prop:normalityim}).

\medskip
 
{\sc Notation and Terminology}

All varieties in this article are defined over the complex 
number field $\mathbb{C}$. 
Let $V$ be a smooth variety. 
We denote by $K_V$ and $\omega_V$ a canonical divisor 
and the dualizing sheaf, respectively, of $V$.  
For a divisor $D$, we denote by $\mathcal{O} (D)$ the coherent 
sheaf associated to $D$.  
For a coherent sheaf $\mathcal{F}$ on $V$, 
we denote by $H^i (\mathcal{F})$, $h^i (\mathcal{F})$, and 
$\chi (\mathcal{F})$,  
the $i$-th cohomology group of $\mathcal{F}$, 
its dimension $\dim_{\mathbb{C}} H^i (\mathcal{F})$, and 
the Euler characteristic $\sum (-1)^i h^i (\mathcal{F})$, 
respectively.  
We denote by $S^n (\mathcal{F})$ and $\bigwedge^n \mathcal{F}$ 
the $n$-th symmetric product and the $n$-th exterior product, 
respectively, of $\mathcal{F}$. 
Let $f: V \to W$ be a morphism to a smooth variety $W$, and 
$D$, a divisor on $W$. We denote by $f^* (D)$   
the total transform of $D$. 

The symbols $\sim$ and $\sim_{\mathrm{num}}$ mean 
the linear equivalence and the numerical equivalence, respectively, 
of two divisors. 
If $D$ and $D^{\prime}$ are two divisors on $V$ and 
$D- D^{\prime}$ is a non-negative divisor, 
we write $D \succeq D^{\prime}$.  

For a smooth algebraic surface $S$,  
we denote by $c_1 (S)$, $p_g(S)$, and $q(S)$, the 
first Chern class, the geometric genus, and the irregularity 
of $S$, respectively.

\section{Some numerical restrictions}  

Let $S$ be a minimal algebraic surface with $c_1^2= 9$ and   
$\chi=5$ whose canonical class is divisible by $3$ in 
the cohomology group $H^2 (S, \mathbb{Z})$. 
We take a divisor $L$ such that $K=K_S \sim 3L$. 
In this section, as a preliminary, we shall find some 
restrictions to numerical invariants 
associated to the divisor $L$. 

Let us begin with the dimension $h^0(\mathcal{O}_S (2L))$.

\begin{lemma}  \label{lemma:ineq2l}
$3 \leq h^0(\mathcal{O}_S (2L)) \leq 5$. 
\end{lemma}

Proof. By the Riemann-Roch theorem, we see that  
\begin{equation} \label{eql:firstinequality}
h^0(\mathcal{O}_S(L)) + h^0(\mathcal{O}_S(2L)) \geq 4.
\end{equation}
By this together with 
$h^0(\mathcal{O}_S(L)) \leq h^0(\mathcal{O}_S(2L))$, 
we obtain $2 \leq h^0(\mathcal{O}_S(2L))$. 
But if $h^0(\mathcal{O}_S(2L)) = 2$, then 
by (\ref{eql:firstinequality})  
we must have $2 \leq h^0(\mathcal{O}_S(L))$, which 
contradicts 
$h^0(\mathcal{O}_S(2L)) \geq 2 h^0(\mathcal{O}_S(L)) -1$. 
Thus we obtain $3 \leq h^0(\mathcal{O}_S(2L))$. 
To obtain the remaining inequality, use 
$h^0(\mathcal{O}_S(6L)) = \chi(\mathcal{O}_S) + K^2 = 14$ 
and $h^0(\mathcal{O}_S(6L)) \geq 3 h^0(\mathcal{O}_S(2L)) -2$. 
\qed 
   
Let $\varPhi_{|2L|} : S - - \to \mathbb{P}^{l_2}$ 
be the rational map associated to the 
linear system $|2L|$, where 
$l_2 = h^0(\mathcal{O}_S(2L)) -1$. 
We have two cases:
the case where $\varPhi_{|2L|}$ is composite with a pencil and 
the case where  $\varPhi_{|2L|}$ is not composite with a pencil.
First, we study the former case.  

\begin{lemma} \label{lm:compwithap}
Assume that the rational map $\varPhi_{|2L|}$ is 
composite with a pencil $\mathcal{P}$. 
Then $h^0(\mathcal{O}_S(L)) = 2$ and $h^0(\mathcal{O}_S(2L)) = 3$ 
hold. Moreover $|L|$ has no fixed component, and 
the pencil $\mathcal{P}$ is given by 
$\varPhi_{|L|}: S - - \to \mathbb{P}^1$. 
\end{lemma}

Proof. Assume that $\varPhi_{|2L|}$ is 
composite with a pencil $\mathcal{P}$. 
Since $S$ is regular and $|2L|$ is complete, 
there exists an effective divisor $D_2$ of $S$ such that 
$h^0(\mathcal{O}_S(D_2)) \geq 2$ and 
$|2L| = |l_2 D_2| + F_2$, where $F_2$ is the fixed part of $|2L|$, 
and $l_2$ is as in the definition of $\varPhi_{|2L|}$.  
Naturally, we have 
\begin{equation} \label{eql:twicelsquared}
2 = 2 L^2 = l_2 D_2 L + F_2 L.
\end{equation}

Assume that we have $D_2 L = 0$. Then we have 
$F_2 L =2$, which together with 
$2LD_2 = l_2 D_2^2 + D_2 F_2$ and 
$2L F_2 = l_2 D_2 F_2 + F_2^2$ implies 
$F_2^2 = 4$ and $D_2^2 = D_2 F_2 = 0$. 
Then by Hodge's Index Theorem, 
we obtain $D_2 = 0$,  
which contradicts the definition of the divisor $D_2$. 

Thus $D_2 L > 0$ holds.  
Since we have $l_2 \geq 2$ by Lemma \ref{lemma:ineq2l},  
we see from this together with (\ref{eql:twicelsquared}) that 
$l_2 = 2$, $D_2 L = 1$, and $F_2 L = 0$.   
In particular, we obtain    
$2 = 2 L D_2 = l_2 D_2^2 + D_2 F_2$. 
But $D_2^2$ is odd, since $D_2K = 3$. 
Thus this implies 
$D_2^2 =1$ and $D_2 F_2 = F_2^2 = 0$, 
hence $F_2 = 0$.  
Thus we obtain $2L \sim l_2 D_2 + F_2 \sim 2 D_2$. 
This however implies $L \sim D_2$, 
since by \cite[Theorem 4]{remarks'''''}
the surface $S$ has no torsion. 
Since $l_2 = h^0(\mathcal{O}_S (2L)) -1$, 
the assertion follows from this linear equivalence and 
$h^0(\mathcal{O}_S (2L)) \geq 2 h^0 (\mathcal{O}_S (L)) -1$.  
\qed

Next, we study the latter case. 
In what follows, we denote by $|M_2|$ and $F_2$ the 
variable part and the fixed part, respectively,  
of the linear system $|2L|$. 
We also denote by $p_2 : \tilde{S}_2 \to S$ the shortest 
composite of quadric transformations such that 
the variable part of $p_2^* |M_2|$ is 
free from base points. 

\begin{lemma}  \label{lm:notcompwap}
Assume that the rational map $\varPhi_{|2L|}$ is not 
composite with a pencil. 
Then $h^0(\mathcal{O}_S(2L)) = 3$, 
$h^0(\mathcal{O}_S(L)) = 1$, and 
$h^1(\mathcal{O}_S(L)) = 0$ hold. 
Moreover, the inequality $2 \leq \tilde{M}_2^2 \leq 4$ holds, 
where $|\tilde{M}_2|$ is the variable part of 
the linear system $p_2^* |M_2|$. 
\end{lemma}

Proof. 
By Lemma \ref{lemma:ineq2l} we have 
$2 \leq l_2 \leq 4$, where 
$l_2 = h^0(\mathcal{O}_S(2L)) - 1$. 
Since we have assumed that $\varPhi_{|2L|}$ is not 
composite with a pencil, a general member $\tilde{M}_2$
of $|\tilde{M}_2|$ is a smooth irreducible curve on $\tilde{S}_2$. 
In what follows, we assume that $\tilde{M}_2$ is general, hence smooth,  
and define the divisors $E_2$ and $\varepsilon_2$ 
by $p_2^* |M_2| = |\tilde{M}_2| + E_2$ and  
$\tilde{K} = K_{\tilde{S}_2} \sim p_2^* (3L) + \varepsilon_2$, 
respectively.  
 
First, let us show that $l_2 \leq 3$. 
By the Serre duality and 
$h^0(\mathcal{O}_{\tilde{S}_2}(p_2^*(L +F_2) + \varepsilon_2 + E_2))
= h^0(\mathcal{O}_S (K- M_2)) $
we have 
$h^2(\mathcal{O}_{\tilde{S}_2}(\tilde{M}_2)) 
= h^0(\mathcal{O}_{\tilde{S}_2}(\tilde{K}-M_2))< p_g(S)= 4.$ 
From this together with the standard exact sequence 
\[
0 \to \mathcal{O}_{\tilde{S}_2} 
   \to \mathcal{O}_{\tilde{S}_2} (\tilde{M}_2)
   \to \mathcal{O}_{\tilde{M}_2} (\tilde{M}_2)
   \to 0,  
\]
we see easily that 
$h^1(\mathcal{O}_{\tilde{M}_2} (\tilde{M}_2)) \geq 1$. 
Thus applying Clifford's theorem for  
$\tilde{M}_2 |_{\tilde{M}_2}$, we obtain
\begin{equation} \label{eql:selfinttwiceL}
2(l_2 - 1) \leq \tilde{M}_2^2 \leq 
\tilde{M}_2^2 + \tilde{M}_2 E_2 + M_2 F_2 + 2L F_2 = (2L)^2 = 4,
\end{equation}
hence in particular $l_2 \leq 3$. 

Assume that we have $l_2 = h^0 (\mathcal{O}_S (2L)) - 1 = 2$.  
Then by the Riemann-Roch theorem,  
we have $h^0 (\mathcal{O}_S (L)) 
= h^1(\mathcal{O}_S (L)) - h^0(\mathcal{O}_S (2L)) + 4 
\geq 1$. 
And also, we have $3 = h^0 (\mathcal{O}_S (2L)) 
\geq 2h^0(\mathcal{O}_S (L)) - 1$, hence 
$2 \geq h^0 (\mathcal{O}_S (L))$. 
The case $h^0 (\mathcal{O}_S (L)) = 2$ however is impossible, since 
we have assumed that $\varPhi_{|2L|}$ is not composite with a pencil. 
Thus we obtain  
$h^0 (\mathcal{O}_S (L)) = 1$ and $h^1 (\mathcal{O}_S (L)) = 0$. 
Moreover, the inequality $2 \leq \tilde{M}_2^2 \leq 4$ follows 
from (\ref{eql:selfinttwiceL}), hence as in the assertion. 
Therefore, we only need to rule out the case $l_2 = 3$. 

So assume that we have $l_2 = 3$. In this case 
we obtain by (\ref{eql:selfinttwiceL}) that
$\tilde{M}_2 E_2 = M_2 F_2  = 2L F_2 =0$,  
which implies the base point freeness of the linear system $|2L|$. 
Since $S$ is of general type, we infer easily from this that  
$\deg \varPhi_{|2L|} = \deg \varPhi_{|2L|} (S) = 2$. 
Thus we have two cases:
\medskip

Case A: the image $\varPhi_{|2L|} (S) \subset \mathbb{P}^3$ is a smooth quadric;

Case B: the image $\varPhi_{|2L|} (S) \subset \mathbb{P}^3$ is a quadric cone. 
\medskip

\noindent
In what follows, we put $g= \varPhi_{|2L|}$. 
We shall rule out the two cases separately. 

{\bf Case A.} 
Assume that $\varPhi_{|2L|} (S)$ is a smooth quadric. 
Then the image $\varPhi_{|2L|} (S)$ is the Hirzebruch surface 
$\varSigma_0$ of degree $0$ embedded by 
$|\varDelta_0 +  \varGamma|$, where 
$\varDelta_0$ and $\varGamma$ denote the minimal section 
and a fiber of the Hirzebruch surface $\varSigma_0$, respectively. 
Let $R$ and $B= g_* (R)$ denote the ramification divisor 
and the branch divisor of the generically 
two-to-one morphism $g: S \to \varSigma_0$, respectively.  
Then since $2L \sim g^*(\varDelta_0 + \varGamma)$, 
we see easily that $R \sim 7L$, hence  
$B \varDelta_0 = B \varGamma = 7$. 
This however is impossible, because $B$ needs to be 
linearly equivalent to twice 
a divisor on $\varSigma_0$. Thus Case A does not occur.  
   
{\bf Case B.}
Assume that $\varPhi_{|2L|} (S)$ is a quadric cone. 
Then the image $\varPhi_{|2L|} (S)$ is the image of 
the morphism 
$\varPhi_{|\varDelta_0 + 2\varGamma|} : \varSigma_2 \to \mathbb{P}^3$,
where $\varSigma_2$ is a Hirzebruch surface of degree $2$, and 
$\varDelta_0$ and $\varGamma$ are its minimal section and a fiber, 
respectively. 
Let $p_2^{\prime}: S_2^{\prime} \to S$ be the shortest 
composite of the quadric transformations such that 
$g \circ p_2^{\prime}$ lifts to a morphism 
$g^{\prime} : S_2^{\prime} \to \varSigma_2$. 
We denote by $K^{\prime} = K_{S_2^{\prime}}$ a canonical divisor 
of $S_2^{\prime}$, and define the divisor $\varepsilon_2^{\prime}$ 
by $K^{\prime} \sim {p_2^{\prime}}^* (3L) + \varepsilon_2^{\prime}$. 
We also denote by $R$ and $B = g^{\prime}_* (R)$ the ramification divisor 
and the branch divisor of the generically two-to-one morphism 
$g^{\prime} : S_2^{\prime} \to \varSigma_2$. 

Since $\varepsilon_2^{\prime}$ is contracted by $g \circ p_2^{\prime}$, 
there exists a natural number $\nu$ such that 
$g^{\prime}_* (\varepsilon_2^{\prime}) = \nu \varDelta_0$. 
Then from 
${p_2^{\prime}}^* (3L) + \varepsilon_2^{\prime} \sim 
{g^{\prime}}^* (-2\varDelta_0  - 4 \varGamma) + R$ and 
$p_2^{\prime} (2L) \sim {g^{\prime}}^* (\varDelta_0 + 2 \varGamma)$ 
we infer that 
$B \varDelta_0 = -2\nu$ and $B \varGamma = 7 + \nu$. 
Since $B$ is linearly equivalent to twice a divisor on $\varSigma_2$, 
this implies $\nu \geq 1$, hence $B \varDelta_0 <0$. 
Thus $\varDelta_0$ is a component of the branch divisor $B$.
In particular, we have $\nu = 1$, from which we see that the 
multiplicity in $\varepsilon_2^{\prime}$ of the $(-1)$-curve
appearing at the last quadric transformation in $p_2^{\prime}$ 
is equal to $1$. Thus   
$p_2^{\prime}: S_2^{\prime} \to S$ is a blowing up at one point, 
and $\varepsilon_2^{\prime}$ is a $(-1)$-curve.  
Then by 
${p_2^{\prime}}^* (2L) \sim {g^{\prime}}^* (\varDelta_0 + 2\varGamma) $
we obtain 
$2 ({p_2^{\prime}}^* L - \varepsilon_2^{\prime} - {g^{\prime}}^*\varGamma)) 
\sim 0$. This implies the linear equivalence 
${p_2^{\prime}}^* L \sim \varepsilon_2^{\prime} + {g^{\prime}}^* \varGamma$, 
since by \cite[Theorem 4]{remarks'''''} our surface $S$ has no torsion. 
Thus we obtain 
$h^0 (\mathcal{O}_S (L)) 
\geq h^0(\mathcal{O}_{\varSigma_2} (\varGamma)) =2$.
This however is impossible, since we have 
$h^0 (\mathcal{O}_S (2L))  = 4$ and 
$4 = h^0 (\mathcal{O}_S (3L))  
\geq h^0 (\mathcal{O}_S (2L))  + h^0 (\mathcal{O}_S (L)) - 1$. 
Thus Case B does not occur.  
This concludes the proof of Lemma \ref{lm:notcompwap}. \qed


\section{Study of the map  $\varPhi_{|2L|}$}

In this section, we shall study the map $\varPhi_{|2L|}$, and  
rule out the case where $\varPhi_{|2L|}$ is 
composite with a pencil. 
Assume that the rational map 
$\varPhi_{|2L|}$ is composite 
with a pencil $\mathcal{P}$. 
Then by Lemma \ref{lm:compwithap}, we have 
$h^0(\mathcal{O}_S(2L)) = 3$ and $h^0(\mathcal{O}_S(L)) = 2$. 
The linear system $|L|$ has a unique base point, which is simple.  
Moreover, 
since $h^0 (\mathcal{O}_L) = 1$ holds, 
$\mathcal{P}$ is a pencil of curves of genus $3$,   
whose members correspond to fibers of  
$\varPhi_{|L|}: S - - \to \mathbb{P}^1$. 
Let $p : \tilde{S} \to S$ be the blow up of $S$ at 
the base point of $|L|$, and $E$, its exceptional curve. 
We denote by 
$f = \varPhi_{|p^*L - E|}: \tilde{S} \to B = \mathbb{P}^1$ 
the morphism associated to the linear system $|p^*L - E|$.  

Since the multiplication map 
$S^3 (H^0(\mathcal{O}_S(L))) \to H^0(\mathcal{O}_S(3L))$ is 
surjective, the canonical map $\varPhi_{|K|}: S - - \to \mathbb{P}^3$ 
is also composite with the pencil $\mathcal{P}$.  
Thus we are in Case (ii) of \cite[Proposition 1.7]{onbicanonicalmaps}. 
In particular, 
any general member of $|L|$ is non-hyperelliptic, 
and all the fibers of 
$f: \tilde{S} \to B$ are $2$-connected. 
Therefor, we can utilize the structure theorem given in 
\cite{fibrationsI'} for $2$-connected non-hyperelliptic fibrations 
of genus $3$.  

In what follows, we put $\tilde{L} = p^* L - E$ and 
$\tilde{K} = K_{\tilde{S}} = p^*(3L) + E$, and denote by 
$\omega_{S|B} = \mathcal{O}_S (\tilde{K} - f^*{K_B})$ 
the relative canonical sheaf of the fibration 
$f: \tilde{S} \to B$.  
Moreover we denote 
by $V_n = f_* (\omega_{\tilde{S}|B}^{\otimes n})$ the direct image by $f$ of 
the sheaf $\omega_{\tilde{S}|B}^{\otimes n}$. 
Recall that for any integer $n \geq 2$ we have 
\[
\mathrm{rk}\, V_n = 4n -2 , \qquad \deg V_n = 7 + 12 n(n-1). 
\]  
The latter equality on $\deg V_n$ is valid also for 
$n= 1$, but for the former equality on $\mathrm{rk}\, V_n$, we have 
instead $\mathrm{rk}\, V_1 = 3$ for $n=1$.  

\begin{lemma} \label{lm:v1v2v4}
The following hold$:$
\smallskip

\noindent
$1$$)$ 
$V_1 \simeq 
\mathcal{O}_B(1)^{\oplus 2} \oplus \mathcal{O}_B(5)$,

\noindent
$2$$)$ 
$V_2 \simeq \left(\bigoplus_{k =2}^{4}\mathcal{O}_B(k) \right) 
\oplus
\mathcal{O}_B(6)^{\oplus 2} \oplus \mathcal{O}_B(10)$,  

\noindent
$3$$)$ 
$V_4 \simeq \left(\bigoplus_{k =4}^{14}\mathcal{O}_B(k) \right) 
\oplus
\mathcal{O}_B(16)^{\oplus 2} \oplus \mathcal{O}_B(20)$.
\end{lemma}

Proof. 
Recall that we have $\mathrm{rk}\, V_1 = 3$ and 
$\deg V_1 =7$. Thus we can put 
$V_1 \simeq \bigoplus_{i =0}^{2} \mathcal{O}_B(a_i)$, 
where $a_0 \leq a_1 \leq a_2$ and $\sum_{i=0}^2 a_i =7$. 
Moreover we have
$\omega_{\tilde{S}|B} 
\simeq \mathcal{O}_{\tilde{S}} (\tilde{K} - f^* K_B) 
\simeq \mathcal{O}_{\tilde{S}} (5 \tilde{L} + 4 E)$. 
Thus we obtain  
\[
h^0(V_1 \otimes \mathcal{O}_B (-k)) = 
h^0(\mathcal{O}_{\tilde{S}} ((5-k)\tilde{L} + 4 E)) = 
h^0(\mathcal{O}_S ((5-k)L)) 
\] 
for any $k \geq 1$, from which we infer 
$h^0 (V_1 \otimes \mathcal{O}_B (-1))) 
- h^0 (V_1 \otimes \mathcal{O}_B (-2))) = 3$. 
This implies $a_i \geq 1$ for all $0 \leq i \leq 2$. 
Since 
$h^0 (V_1 \otimes \mathcal{O}_B (-k))) 
- h^0 (V_1 \otimes \mathcal{O}_B (-(k+1)))$ 
is equal to the numbers of $i$'s satisfying  
$a_i \geq k$, using Lemma \ref{lm:compwithap}, 
we obtain the assertion $1$). 

The assertions $2$) and $3$) can be proved exactly in the same way. 
For these two, use  
$\omega_{\tilde{S}|B}^{\otimes 2} 
\simeq \mathcal{O}_{\tilde{S}} (10 \tilde{L} + 8 E)$ and 
$\omega_{\tilde{S}|B}^{\otimes 4} 
\simeq \mathcal{O}_{\tilde{S}} (20 \tilde{L} + 16 E)$. \qed

In what follows, we denote by 
$X_0$, $X_1$, and $X_2$ local bases of 
the direct summands $\mathcal{O}_B(1)$, $\mathcal{O}_B(1)$, 
and $\mathcal{O}_B(5)$, respectively,  of the sheaf $V_1$. 
We also denote by 
$S_0$, $S_1$, $S_2$, $T_0$, $T_1$, and $U_0$ 
local bases of the direct summands 
$\mathcal{O}_B(2)$, $\mathcal{O}_B(3)$, $\mathcal{O}_B(4)$, 
$\mathcal{O}_B(6)$, $\mathcal{O}_B(6)$, and 
$\mathcal{O}_B(10)$, respectively, 
of the sheaf $V_2$. 
By Lemma \ref{lm:v1v2v4} we have
\[
S^2(V_1) \simeq 
\mathcal{O}_B(2)^{\oplus 3} \oplus
\mathcal{O}_B(6)^{\oplus 2} \oplus
\mathcal{O}_B(10),  
\]
where the local bases of the direct summands are given by 
$X_0^2$, $X_0 X_1$, $X_1^2$, $X_0 X_2$, $X_1 X_2$, and $X_2^2$, 
respectively.   
With these local bases, the multiplication morphism 
$\sigma_2 : S^2 (V_1) \to V_2$ is expressed by a  
$6 \times 6$ matrix $A$ in the following form: 
\begin{equation}    \label{eql:AandAprime}
A = 
\begin{pmatrix}
A^{\prime} & O_3 \\
* & I_3
\end{pmatrix}, 
\quad 
\text{where} 
\quad
A^{\prime} = 
\begin{pmatrix}
a_0      & a_1      & a_2 \\
\alpha_0 & \alpha_1 & \alpha_2 \\
\beta_0  & \beta_1  & \beta_2  
\end{pmatrix}.  
\end{equation}
Here $O_3$ and $I_3$ denote the $3 \times 3$ zero matrix 
and the $3 \times 3$ identity matrix, respectively, and  
$a_i \in H^0(\mathcal{O}_B )$,   
$\alpha_j \in H^0(\mathcal{O}_B (1))$, and  
$\beta_k \in H^0(\mathcal{O}_B (2))$ 
are global sections for each $0 \leq i,\, j,\, k \leq 2$.

Let us describe the $5$-tuple for our genus $3$ fibration 
$f: \tilde{S} \to B$. 
For the notion of the $5$-tuple, see \cite{fibrationsI'}. 
Let $\tau$ be the effective divisor of degree $\deg \tau = 3$ 
on $B$ determined by the 
short exact sequence
\begin{equation} \label{eql:extensionclass}
0 \to S^2 (V_1) \to V_2 \to \mathcal{O}_{\tau} \to 0. 
\end{equation}
Let  
$\mathcal{C} : S^2 (\bigwedge^2 V_1) \to S^2(S^2 (V_1))$ 
be the morphism given by 
$(a \wedge b) (c \wedge d) \mapsto (ac)(bd) - (ad)(bc)$. 
Then the morphism 
$S^2 (\sigma_2) \circ \mathcal{C} : S^2 (\bigwedge^2 V_1) \to S^2 (V_2)$ 
has a locally free cokernel of rank $15$, which we shall denote 
by $\tilde{V}_4 = \mathrm{Cok}\, (S^2 (\sigma_2) \circ \mathcal{C})$.  
We denote by $\mathcal{L}_4^{\prime}$ and $\mathcal{L}_4$ 
the kernel of the natural surjection $\tilde{V_4} \to V_4$ and 
that of the natural morphism $S^4 (V_1) \to V_4$, respectively. 
Then we obtain the natural inclusion morphism 
\begin{equation} \label{eql:l4prime}
\mathcal{L}_4^{\prime} 
\simeq (\det V_1 ) \otimes \mathcal{O}_B (\tau) 
\simeq \mathcal{O}_B (10)
\to \tilde{V}_4 . 
\end{equation}
With the notation above, $B$, $V_1$, $\tau$, 
(\ref{eql:extensionclass}), and (\ref{eql:l4prime}) 
form the admissible $5$-tuple associated to 
our fibration $f: \tilde{S} \to B$.    

By Lemma \ref{lm:v1v2v4} we have 
$\bigwedge^2 V_1 \simeq 
\mathcal{O}_B(2) \oplus \mathcal{O}_B(6)^{\oplus 2}$ 
and 
$S^2(\bigwedge^2 V_1) \simeq 
\mathcal{O}_B(4)  
\oplus \mathcal{O}_B(8)^{\oplus 2} \oplus \mathcal{O}_B(12)^{\oplus 3}$. 
We decompose each of the five sheaves 
$S^2(\bigwedge^2 V_1)$, $S^2(V_1)$, 
$S^2 (S^2 (V_1))$, $V_2$, and $S^2(V_2)$ 
into the lower degree part $(\mathrm{L})$ 
and the higher degree part $(\mathrm{H})$ 
as follows: 
\begin{align}
S^2 (\bigwedge^2 V_1)  
&= \left[ \mathcal{O}_B(4) \right] \oplus 
\left[ \mathcal{O}_B(8)^{\oplus 2} \oplus \mathcal{O}_B(12)^{\oplus 3} \right]  
\notag \\
&= S^2 (\bigwedge^2 V_1)^{(\mathrm{L} )} \oplus S^2 (\bigwedge^2 V_1)^{(\mathrm{H})},
\notag \\ 
S^2(V_1) &= \left[ \mathcal{O}_B(2)^{\oplus 3} \right] \oplus 
\left[ \mathcal{O}_B(6)^{\oplus 2} \oplus \mathcal{O}_B(10) \right] \notag \\ 
&= S^2(V_1)^{(\mathrm{L} )} \oplus S^2(V_1)^{(\mathrm{H} )} ,
\notag \\ 
S^2(S^2 (V_1)) 
&= 
\left[ S^2 (S^2(V_1)^{(\mathrm{L} )}) \right] \oplus
\left[ 
\left(
S^2(V_1)^{(\mathrm{L} )} \otimes S^2(V_1)^{(\mathrm{H} )}
\right) 
\oplus S^2(S^2(V_1)^{(\mathrm{H} )}) 
\right] 
\notag \\
&= 
S^2(S^2 ( V_1))^{(\mathrm{L} )}   
\oplus S^2(S^2( V_1))^{(\mathrm{H})} ,
\notag \\
V_2 &=  \left[ \bigoplus_{k=2}^{4} \mathcal{O}_B(k) \right] \oplus 
\left[ \mathcal{O}_B(6)^{\oplus 2} \oplus \mathcal{O}_B(10) \right] \notag \\
&= V_2^{(\mathrm{L} )} \oplus V_2^{(\mathrm{H} )} ,
\notag \\
S^2 (V_2)  &=
\left[ S^2 (V_2^{(\mathrm{L} )}) \right]
\oplus
\left[ \left( V_2^{(\mathrm{L} )} \otimes V_2^{(\mathrm{H} )} \right) 
\oplus
S^2 (V_2^{(\mathrm{H} )}) 
\right] 
\notag \\
&=  S^2(V_2)^{(\mathrm{L} )} \oplus S^2(V_2)^{(\mathrm{H} )} ,
\notag 
\end{align}
where in each expression the first $[ \quad ]$ term corresponds 
to the lower degree part $(\mathrm{L})$, and 
the second $[ \quad  ]$ term corresponds 
to the higher degree part $(\mathrm{H})$.   

Let 
$\gamma 
: S^2 (\bigwedge^2 V_1)^{(\mathrm{L} )} 
\simeq \mathcal{O}_B (4) \to S^2(V_2)^{(\mathrm{L} )}$ be 
the composition of the morphism 
$\mathcal{C} |_{S^2 (\bigwedge^2 V_1)^{(\mathrm{L} )}}
: S^2 (\bigwedge^2 V_1)^{(\mathrm{L} )} \to 
S^2 ( S^2 (V_1))^{(\mathrm{L} )} $ and 
the morphism 
$S^2(A^{\prime}) : 
S^2(S^2 ( V_1))^{(\mathrm{L} )} = S^2 (S^2(V_1)^{(\mathrm{L} )}) 
\to 
S^2(V_2)^{(\mathrm{L} )} = S^2 (V_2^{(\mathrm{L} )}) 
$, 
where $A^{\prime}$ is the $3 \times 3$ matrix given in 
(\ref{eql:AandAprime}). 

\begin{lemma}  \label{lm:homl4cok}
$\mathrm{Hom}_{\mathcal{O}_B}(\mathcal{L}_4^{\prime}, \mathrm{Cok}\, \gamma ) 
\neq \{ 0 \}$. 
\end{lemma}

Proof. 
Note that  by (\ref{eql:AandAprime}) we have  
$(S^2 (\sigma_2) \circ \mathcal{C}) (S^2 (\bigwedge^2 V_1)^{(\mathrm{H} )})
\subset S^2(V_2)^{(\mathrm{H})}$. 
Thus 
$(S^2 (\sigma_2) \circ \mathcal{C}) : 
S^2 (\bigwedge^2 V_1) \to S^2(V_2)$ induces 
a morphism of $\mathcal{O}_B$-modules 
\[
\gamma^{\prime} : 
\frac{S^2 (\bigwedge^2 V_1)}{S^2 (\bigwedge^2 V_1)^{(\mathrm{H} )}} 
\simeq S^2 (\bigwedge^2 V_1)^{(\mathrm{L} )} 
\to 
\frac{S^2(V_2)}{S^2(V_2)^{(\mathrm{H} )}}
\simeq
S^2(V_2)^{(\mathrm{L} )}. 
\] 
Our morphism $\gamma$ coincides with this $\gamma^{\prime}$, 
when we view $\gamma^{\prime}$ as a morphism 
from $S^2 (\bigwedge^2 V_1)^{(\mathrm{L} )}$ 
to $S^2(V_2)^{(\mathrm{L} )}$.  
Thus by the commutative diagram
\begin{equation}   \label{dgm:admsv4tilde}
\begin{CD}  
0     @>\text{$$}>>
S^2 (\bigwedge^2 V_1)^{(\mathrm{H} )}  @>\text{$$}>>  
S^2(V_2)^{(\mathrm{H})}        @>\text{$$}>>
\frac{S^2(V_2)^{(\mathrm{H})}}{S^2 (\bigwedge^2 V_1)^{(\mathrm{H} )}} 
@>\text{$$}>> 0 \\
@. 
@VV\text{$$}V  
@VV\text{$$}V 
@VV\text{$$}V
@.              \\
0     @>\text{$$}>>
S^2 (\bigwedge^2 V_1) @>\text{$S^2 (\sigma_2) \circ \mathcal{C}$}>>  
S^2(V_2)   @>\text{$$}>> 
\tilde{V}_4  @>\text{$$}>>
0 
\end{CD}
\end{equation}
and $3 \times 3$ Lemma, we obtain the following two 
short exact sequences: 
\begin{gather}
0 
\to 
\frac{S^2 (\bigwedge^2 V_1)}{S^2 (\bigwedge^2 V_1)^{(\mathrm{H} )}} 
\to
\frac{S^2(V_2)}{S^2(V_2)^{(\mathrm{H})}}
\to
\mathrm{Cok}\, \gamma^{\prime} \simeq \mathrm{Cok}\, \gamma
\to
0,    \notag \\
0
\to
\frac{S^2(V_2)^{(\mathrm{H})}}{S^2 (\bigwedge^2 V_1)^{(\mathrm{H} )}} 
\to
\tilde{V}_4
\to
\mathrm{Cok}\, \gamma^{\prime} \simeq \mathrm{Cok}\, \gamma
\to
0.     \label{eql:exctseqv4tilde} 
\end{gather}

Now, assume that we have 
$\mathrm{Hom}_{\mathcal{O}_B}(\mathcal{L}_4^{\prime}, \mathrm{Cok}\, \gamma ) 
= \{ 0 \}$. 
Then by the short exact sequence (\ref{eql:exctseqv4tilde}) above,  
we obtain the surjectivity of the morphism 
$\mathrm{Hom}_{\mathcal{O}_B} 
(
\mathcal{L}_4^{\prime} ,\, 
\frac{S^2(V_2)^{(\mathrm{H})}}{S^2 (\bigwedge^2 V_1)^{(\mathrm{H} )}} 
) 
\to
\mathrm{Hom}_{\mathcal{O}_B} (\mathcal{L}_4^{\prime} ,\, \tilde{V}_4)$. 
This implies that the morphism  
$\frac{S^2(V_2)^{(\mathrm{H})}}{S^2 (\bigwedge^2 V_1)^{(\mathrm{H} )}}
\to \tilde{V}_4$ in (\ref{dgm:admsv4tilde}) 
factors through the inclusion morphism (\ref{eql:l4prime}).  
On the other hand, however, since 
$\sigma_2 |_{S_2(V_1)^{(\mathrm{H})}} : 
S^2(V_1)^{(\mathrm{H})} \to V_2^{(\mathrm{H})}$ is an isomorphism
by (\ref{eql:AandAprime}), 
we have also the surjectivity of the morphism
$V_2 \otimes S^2(V_1)^{(\mathrm{H})} \to 
S^2(V_2)^{(\mathrm{H})} = V_2 \cdot V_2^{(\mathrm{H})}$.   
Then with the help of the commutative diagram (\ref{dgm:admsv4tilde}), 
we find immediately a contradiction to the definition of 
an admissible $5$-tuple 
(see \cite[Condition (iv), Definition 7.10]{fibrationsI'}). 
Thus 
$\mathrm{Hom}_{\mathcal{O}_B}(\mathcal{L}_4^{\prime}, \mathrm{Cok}\, \gamma ) 
= \{ 0 \}$
is impossible. \qed

Note that by Lemma \ref{lm:v1v2v4} we have 
\[
S^2(V_2)^{(\mathrm{L})} \simeq 
\mathcal{O}_B(4) \oplus \mathcal{O}_B(5) \oplus \mathcal{O}_B(6)^{\oplus 2} 
\oplus \mathcal{O}_B(7) \oplus \mathcal{O}_B(8). 
\]
Local bases of the direct summands are given by 
$S_0^2$, $S_0 S_1$, $S_1^2$, $S_0 S_2$, $S_1 S_2$, 
and $S_2^2$, respectively.  
In what follows, we shall compute the sheaf  
$\mathrm{Cok}\, \gamma$, and rule out the case where 
$\varPhi_{|2L|}$ is composite with a pencil. 
For this we divide our argument into several cases, 
normalizing the matrix $A^{\prime}$. 

First, by replacing the bases $X_0$ and $X_1$ of the sheaf $V_1$, we may 
assume $a_1 =1$. 
Then by replacing the bases $S_0$, $S_1$, and $S_2$ of the sheaf $V_2$, 
we may assume $\alpha_1 =0$ and $\beta_1 =0$. 
Then we obtain
\[
A^{\prime} = 
\begin{pmatrix}
a_0      & 1 & a_2 \\
\alpha_0 & 0 & \alpha_2 \\
\beta_0  & 0 & \beta_2  
\end{pmatrix}.  
\]  
We have two cases:
\medskip

Case $1$: $a_0 a_2 \neq 1$; 

Case $2$: $a_0 a_2 = 1$.

\begin{lemma}   \label{lm:ruleout1}
Case $1$ does not occur. 
\end{lemma}

Proof. The composite of   
the morphism $\gamma$ and the natural projection 
$S^2 (V_2)^{(\mathrm{L})} \to \mathcal{O}_B(4)$ coincides with 
$(a_0 a_2 - 1) \times: 
S^2(\bigwedge^2 V_1)^{(\mathrm{L})} \simeq \mathcal{O}_B(4) 
\to \mathcal{O}_B(4)$. 
Thus if we are in Case $1$, then the image $\mathrm{Im}\, \gamma$ 
is a direct summand of $S^2 (V_2)^{(\mathrm{L})}$. Thus we obtain 
\[
\mathrm{Cok}\, \gamma \simeq  
\mathcal{O}_B(5) \oplus \mathcal{O}_B(6)^{\oplus 2} 
\oplus \mathcal{O}_B(7) \oplus \mathcal{O}_B(8),  
\] 
which contradicts Lemma \ref{lm:homl4cok}. \qed 

Let us study Case $2$. 
In this case, the composite of   
the morphism $\gamma$ and the natural projection 
$S^2 (V_2)^{(\mathrm{L})} \to \mathcal{O}_B (4)$ is a zero morphism. 
Thus $\gamma$ is a composite of a morphism    
\[
\gamma_0: S^2(\bigwedge^2 V_1)^{(\mathrm{L})}
\simeq \mathcal{O}_B(4) \to 
\mathcal{F}_0 =  
\mathcal{O}_B(5) \oplus \mathcal{O}_B(6)^{\oplus 2} 
\oplus \mathcal{O}_B(7) \oplus \mathcal{O}_B(8)
\]
and the natural inclusion 
$\mathcal{F}_0 \to S^2 (V_2)^{(\mathrm{L})}$, and 
we find 
$\mathrm{Cok}\, \gamma 
\simeq \mathcal{O}_B(4) \oplus \mathrm{Cok}\, \gamma_0$. 
By replacing the bases $X_0$ and $X_1$ by their multiples by 
non-zero constants, we may assume $a_0 = a_2 = 1$. 
Then by the short exact sequence 
\begin{equation}       \label{eql:exactf0gm0}
0 \to S^2(\bigwedge^2 V_1)^{(\mathrm{L})}  
\to \mathcal{F}_0 \to \mathrm{Cok}\, \gamma_0 \to 0, 
\end{equation}
we obtain
\begin{align}
h^0 (\mathrm{Cok}\, \gamma_0 \otimes \mathcal{O}_B (-6)) &= 8 \notag \\
h^0 (\mathrm{Cok}\, \gamma_0 \otimes \mathcal{O}_B (-7)) &=  
3 + \dim \mathrm{Ker}\, 
((\alpha_0 + \alpha_2) \times ),      \label{eql:dimcok-7}
\end{align}
where 
$(\alpha_0 + \alpha_2) \times :
H^1 (\mathcal{O}_B (-3)) \to H^1 (\mathcal{O}_B (-2))$ is 
the morphism induced by the multiplication 
morphism by $\alpha_0 + \alpha_2$ of sheaves.  

Case $2$ splits into two cases:
\medskip

Case $2$--$1$: $\alpha_0 + \alpha_2 \neq 0 \in H^0 (\mathcal{O}_B (1))$; 

Case $2$--$2$: $\alpha_0 + \alpha_2 = 0 \in H^0 (\mathcal{O}_B (1))$.
\medskip

{\bf Case $2$--$1$}. Let us study Case $2$--$1$. 
In this case the morphism 
$(\alpha_0 + \alpha_2) \times : 
H^1 (\mathcal{O}_B (-3)) \to H^1 (\mathcal{O}_B (-2))$ 
is surjective. Therefor 
by (\ref{eql:exactf0gm0}) and (\ref{eql:dimcok-7}) we obtain 
\begin{align}
h^0 (\mathrm{Cok}\, \gamma_0 \otimes \mathcal{O}_B (-7)) &= 4 \notag \\
h^0 (\mathrm{Cok}\, \gamma_0 \otimes \mathcal{O}_B (-8)) &=  
1 + \dim \mathrm{Ker}\, 
( \,^t (\alpha_0 + \alpha_2 , \, 
\alpha_0  \alpha_2 , \,
\beta_0 + \beta_2) \times),   \label{eql:cokgamma0-8}
\end{align}
where 
$ \,^t (\alpha_0 + \alpha_2 , \, 
\alpha_0  \alpha_2 , \,
\beta_0 + \beta_2) \times : 
H^1 (\mathcal{O}_B (-4)) \to 
H^1 (\mathcal{O}_B (-3) \oplus \mathcal{O}_B (-2)^{\oplus 2})$ 
is the morphism induced by the multiplication morphism by 
$\,^t (\alpha_0 + \alpha_2 , \, 
\alpha_0  \alpha_2 , \,
\beta_0 + \beta_2)$ of sheaves. 

Case $2$--$1$ splits into two cases:
\medskip

Case $2$--$1$--$1$: 
$\alpha_0 + \alpha_2$, $\alpha_0  \alpha_2$, and 
$\beta_0 + \beta_2$ have no common zero; 

Case $2$--$1$--$2$: 
$\alpha_0 + \alpha_2$, $\alpha_0  \alpha_2$, and 
$\beta_0 + \beta_2$ have a common zero.  

\begin{lemma}    \label{lm:ruleout2-1-1}
Case $2$--$1$--$1$ does not occur. 
\end{lemma}

Proof. 
Assume that we are in Case $2$--$1$--$1$. 
Let us denote by 
$\gamma_0^{(1)} : S^2 (\bigwedge^2 V_1)^{(\mathrm{L})} \simeq 
\mathcal{O}_B(4) \to \mathcal{O}_B(5) \oplus  \mathcal{O}_B(6)^{\oplus 2}$
the multiplication morphism by 
$\,^t (\alpha_0 + \alpha_2 , \, 
\alpha_0  \alpha_2 , \,
\beta_0 + \beta_2)$. 
Then both $\mathrm{Cok}\, \gamma_0$ and 
$\mathrm{Cok}\, \gamma_0^{(1)}$ are locally free, and we have 
$\mathrm{rk}\, \mathrm{Cok}\, \gamma_0 = 4$ and 
$\mathrm{rk}\, \mathrm{Cok}\, \gamma_0^{(1)} = 2$. 
Put 
$\mathrm{Cok}\, \gamma_0 \simeq 
\bigoplus_{i=0}^3 \mathcal{O}_B (b_i)$ where 
$b_0 \leq b_1 \leq b_2$. 
Then by the short exact sequence
\[
0 \to 
S^2 (\bigwedge^2 V_1)^{(\mathrm{L})} \to 
\mathcal{O}_B(5) \oplus \mathcal{O}_B(6)^{\oplus 2} \to 
\mathrm{Cok}\, \gamma_0^{(1)} \to 0 , 
\]
we obtain 
$h^0 (\mathrm{Cok}\, \gamma_0^{(1)} \otimes \mathcal{O}_B (-6)) = 3$ 
and 
$h^0 (\mathrm{Cok}\, \gamma_0^{(1)} \otimes \mathcal{O}_B (-7)) = 1$. 
From these together with 
$\deg \mathrm{Cok}\, \gamma_0^{(1)} = 13$, we infer 
$\mathrm{Cok}\, \gamma_0^{(1)} \simeq 
\mathcal{O}_B (6) \oplus \mathcal{O}_B (7)$, which 
in tern together with (\ref{eql:cokgamma0-8}) implies 
$h^0 (\mathrm{Cok}\, \gamma_0 \otimes \mathcal{O}_B (-8))) = 1$. 
This together with (\ref{eql:dimcok-7}) and (\ref{eql:cokgamma0-8}) 
implies $b_0 = 6$ and $b_i \geq 7$ for all $1 \leq i \leq 3$. 
Then since $\deg \mathrm{Cok}\, \gamma_0 =28$, we obtain 
\[
\mathrm{Cok}\, \gamma \simeq 
\mathcal{O}_B(4) \oplus \mathrm{Cok}\, \gamma_0 \simeq 
\mathcal{O}_B(4) \oplus \mathcal{O}_B(6) \oplus 
\mathcal{O}_B(7)^{\oplus 2} \oplus \mathcal{O}_B(8), 
\]
which contradicts Lemma \ref{lm:homl4cok}.  \qed    

\begin{lemma}   \label{lm:ruleout2-1-2}
Case $2$--$1$--$2$ does not occur. 
\end{lemma}

Proof. 
Assume that we are in Case $2$--$1$--$2$. 
Without loss of generality we may assume 
$\alpha_0 \neq 0 \in H^0 (\mathcal{O}_B (1))$. 
Since the three sections 
$\alpha_0 + \alpha_2$, $\alpha_0  \alpha_2$,
and $\beta_0 + \beta_2$ have a common zero, 
there exist a number  $a \in \mathbb{C}$ and 
a section $\lambda \in H^0 (\mathcal{O}_B (1))$ such 
that $\alpha_0 + \alpha_2 = (1+a) \alpha_0$ and 
$\beta_0 + \beta_2  = \lambda \alpha_0 $ hold.  
Note that we have $1+a \neq 0$, since we are in Case $2$--$1$. 
Since $\beta_0 \beta_2 = \beta_0 (\lambda \alpha_0 - \beta_0)$, 
if the two sections $\alpha_0 + \alpha_2$ and $\beta_0 \beta_2$ have 
a common zero $P$, then at this point $P$, the rank of 
$\sigma_2 \otimes k(P)$ drops at least by $2$, which is impossible 
(see \cite{fibrationsI'}).  Thus $\mathrm{Cok}\, \gamma_0$ is locally free. 

Put $\mathrm{Cok}\, \gamma_0 \simeq \bigoplus_{i=0}^3 \mathcal{O}_B (b_i)$, 
where $b_0 \leq b_1 \leq b_2 \leq b_3$.  
Let us denote by 
$\gamma_0^{(1)}: S^2 (\bigwedge^2 V_1)^{(\mathrm{L})} \simeq 
\mathcal{O}_B(4) 
\to \mathcal{O}_B(5) \oplus \mathcal{O}_B (6)^{\oplus 2}$ 
the multiplication morphism by 
$\,^t (\alpha_0 + \alpha_2 ,\, \alpha_0 \alpha_2 ,\, \beta_0 + \beta_2)
= \,^t ( (1+a) \alpha_0 ,\, a\alpha_0^2 ,\, \lambda \alpha_0)$. 
Then denoting by 
$\gamma_0^{(2)}: \mathcal{O}_B(5) 
\to \mathcal{O}_B(5) \oplus \mathcal{O}_B (6)^{\oplus 2}$ 
the multiplication morphism by 
$\,^t ( (1+a),\, a\alpha_0 ,\, \lambda )$, 
we have 
$\gamma_0^{(1)} = \gamma_0^{(2)} \circ (\alpha_0 \times) $, 
where 
$\alpha_0 \times : S^2 (\bigwedge^2 V_1)^{(\mathrm{L})} 
\to \mathcal{O}_B(5)$ is the multiplication morphism 
by $\alpha_0$.   
From this we see easily that the morphism 
$\gamma_0^{(1)} \otimes \mathcal{O}_B(-8): 
S^2 (\bigwedge^2 V_1)^{(\mathrm{L})} \otimes \mathcal{O}_B(-8) 
\to \mathcal{O}_B(-3) \oplus \mathcal{O}_B (-2)^{\oplus 2}$ 
induces a morphism 
$
H^1 (S^2 (\bigwedge^2 V_1)^{(\mathrm{L})} \otimes \mathcal{O}_B(-8)) 
\to H^1 (\mathcal{O}_B(-3) \oplus \mathcal{O}_B (-2)^{\oplus 2})$ 
of rank $4$.  
Thus we obtain 
$h^0 (\mathrm{Cok}\, \gamma_0 \otimes \mathcal{O}_B (-8) ) =2$, 
which together with 
(\ref{eql:dimcok-7}) and (\ref{eql:cokgamma0-8}) implies 
$b_0 = b_1 =6$. 
Then since $\deg \mathrm{Cok} \, \gamma_0 = 28$, we obtain 
\[
\mathrm{Cok}\, \gamma \simeq 
\mathcal{O}_B(4) \oplus \mathrm{Cok}\, \gamma_0 \simeq 
\mathcal{O}_B(4) \oplus \mathcal{O}_B(6)^2 \oplus 
\mathcal{O}_B(b_2) \oplus \mathcal{O}_B(b_3),  
\]
where $(b_2, b_3) = (8, 8)$ or $(7, 9)$, 
which contradicts Lemma \ref{lm:homl4cok}. \qed 
 
{\bf Case $2$--$2$}. Let us study Case $2$--$2$. 
In this case we have 
$\alpha_2 = - \alpha_0 \neq 0 \in H^0 (\mathcal{O}_B (1))$.  
Moreover, $\gamma_0$ is a composite of a morphism    
\[
\gamma_1: S^2(\bigwedge^2 V_1)^{(\mathrm{L})}
\simeq \mathcal{O}_B(4) \to 
\mathcal{F}_1 = 
\mathcal{O}_B(6)^{\oplus 2} 
\oplus \mathcal{O}_B(7) \oplus \mathcal{O}_B(8)
\]
and the natural inclusion 
$\mathcal{F}_1 \to \mathcal{F}_0$, and 
we find 
$\mathrm{Cok}\, \gamma_0  
\simeq \mathcal{O}_B(5) \oplus \mathrm{Cok}\, \gamma_1$. 

Case $2$--$2$ splits into two cases: 
\medskip

Case  $2$--$2$--$1$: $\alpha_0 \alpha_2 = - \alpha_0^2$ and 
$\beta_0 + \beta_2$ have no common zero;

Case  $2$--$2$--$2$: $\alpha_0 \alpha_2 = - \alpha_0^2$ and 
$\beta_0 + \beta_2$ have a common zero.

\begin{lemma}  \label{lm:ruleout2-2-1}
Case  $2$--$2$--$1$ does not occur. 
\end{lemma}

Proof. 
Assume that we are in Case $2$--$2$--$1$. 
Then the sheaf $\mathrm{Cok}\, \gamma_1$ is 
locally free of rank $3$. 
Put $\mathrm{Cok}\, \gamma_1 \simeq \bigoplus_{i=1}^3 \mathcal{O}_B (b_i)$, 
where $b_1 \leq b_2 \leq b_3$.  
Then since the multiplication morphism 
$\,^t (\alpha_0 \alpha_2, \, \beta_0 + \beta_0) \times: 
S^2 (\bigwedge^2 V_1)^{(\mathrm{L})} \to \mathcal{O}_B (6)^{\oplus 2}$  
by $\,^t (\alpha_0 \alpha_2, \, \beta_0 + \beta_0)$ has 
a cokernel isomorphic to $\mathcal{O}_B (8)$, we obtain 
by the short exact sequence 
\[
0 \to S^2(\bigwedge^2 V_1)^{(\mathrm{L})}  
\to \mathcal{F}_1 \to \mathrm{Cok}\, \gamma_1 \to 0  
\]
that 
$h^0 (\mathrm{Cok}\, \gamma_1 \otimes \mathcal{O}_B(-7)) = 5$ and 
$h^0 (\mathrm{Cok}\, \gamma_1 \otimes \mathcal{O}_B(-8)) = 2$, 
which imply 
$b_i \geq 7$ for all $1 \leq i \leq 3$. 
Since $\deg \mathrm{Cok}\, \gamma_1 = 23$, we obtain 
\[
\mathrm{Cok}\, \gamma \simeq 
\mathcal{O}_B(4) \oplus \mathcal{O}_B(5) \oplus 
\mathcal{O}_B(b_1) \oplus
\mathcal{O}_B(b_2) \oplus \mathcal{O}_B(b_3),  
\]
where $(b_1, b_2, b_3) = (7, 7, 9)$ or $(7, 8, 8)$,  
which contradicts Lemma \ref{lm:homl4cok}. \qed 

Let us study Case $2$--$2$--$2$. 
In this case there exists a section 
$\lambda \in H^0(\mathcal{O}_B (1))$ such that 
$\beta_0 + \beta_2 = \lambda \alpha_0$. 
Then since $\beta_0 \beta_2 = \beta_0 (\lambda \alpha_0 - \beta_0)$ holds, 
if the two sections $-\alpha_0^2$ and $\beta_0 \beta_2$ 
have a common zero $P$, 
then at this point $P$, the rank of $\sigma_2 \otimes k(P)$ drops 
at least by $2$, which is impossible.  
Thus $\mathrm{Cok}\, \gamma_1$ is locally free of rank $3$. 
Put $\mathrm{Cok}\, \gamma_1 \simeq \bigoplus_{i=1}^3 \mathcal{O}_B (b_i)$, 
where $b_1 \leq b_2 \leq b_3$. Then by the same short exact sequence 
as in the proof of Lemma \ref{lm:ruleout2-2-1}, we obtain 
\begin{equation}  \label{eql:2-2-2dimcok-6}
h^0 (\mathrm{Cok}\, \gamma_1 \otimes \mathcal{O}_B(-6)) = 8, 
\qquad
h^0 (\mathrm{Cok}\, \gamma_1 \otimes \mathcal{O}_B(-7)) = 5. 
\end{equation}

Case $2$--$2$--$2$ splits into two cases: 
\medskip

Case  $2$--$2$--$2$--$1$: $\alpha_0$ and $\lambda$ have no common zero;

Case  $2$--$2$--$2$--$2$: $\alpha_0$ and $\lambda$ have a common zero.

\begin{lemma}    \label{lm:ruleout2-2-2-1}
Case $2$--$2$--$2$--$1$ does not occur. 
\end{lemma}

Proof. 
Assume that we are in Case $2$--$2$--$2$--$1$. 
Then since the multiplication morphism 
$\,^t (\alpha_0 \alpha_2,\, \beta_0 + \beta_2) \times : 
S^2(\bigwedge^2 V_1)^{(\mathrm{L})}  \to \mathcal{O}_B (6)^{\oplus 2}$ 
is the composite of the two morphisms $\alpha_0 \times : 
S^2(\bigwedge^2 V_1)^{(\mathrm{L})}  \to \mathcal{O}_B (5)$ 
and 
$\,^t (-\alpha_0,\, \lambda) \times : 
\mathcal{O}_B (5) \to \mathcal{O}_B (6)^{\oplus 2}$, 
we see by the same short exact sequence as in the proof of 
Lemma \ref{lm:ruleout2-2-1} that 
$h^0 (\mathrm{Cok}\, \gamma_1 \otimes \mathcal{O}_B(-8)) = 2$, 
which together with (\ref{eql:2-2-2dimcok-6}) 
implies $b_i \geq 7$ for all $1 \leq i \leq 3$. 
Since $\deg \mathrm{Cok}\, \gamma_1 = 23$, we obtain 
\[
\mathrm{Cok}\, \gamma \simeq 
\mathcal{O}_B(4) \oplus \mathcal{O}_B(5) \oplus 
\mathcal{O}_B(b_1) \oplus
\mathcal{O}_B(b_2) \oplus \mathcal{O}_B(b_3),  
\]
where $(b_1, b_2, b_3) = (7, 7, 9)$ or $(7, 8, 8)$,  
which contradicts Lemma \ref{lm:homl4cok}. \qed 

\begin{lemma}   \label{lm:ruleout2-2-2-2}
Case $2$--$2$--$2$--$2$ does not occur.
\end{lemma}

Proof. Assume that we are in Case $2$--$2$--$2$--$2$. 
Then there exists a number $c \in \mathbb{C}$ such that 
$\lambda = c \alpha_0$. 
If we have $c=0$, then we obtain $\beta_2 = -\beta_0$. 
This however is impossible since $\sigma_2 \otimes k(P)$ 
needs to have rank $6$ at a general point $P$ of $B$. 
Thus we obtain $c \neq 0$. 
Note that the multiplication morphism 
$\,^t (-\alpha_0^2, \, \beta_0 + \beta_2) \times : 
S^2(\bigwedge^2 V_1)^{(\mathrm{L})}  \to \mathcal{O}_B (6)^{\oplus 2}$ 
is the composite of the two morphisms $\alpha_0^2 \times : 
S^2(\bigwedge^2 V_1)^{(\mathrm{L})}  \to \mathcal{O}_B (6)$ 
and 
$\,^t (-1,\, c) \times : 
\mathcal{O}_B (6) \to \mathcal{O}_B (6)^{\oplus 2}$. 
Since we have 
$\mathrm{Cok}\, (\,^t (-1,\, c) \times) \simeq \mathcal{O}_B (6)$,  
we see by the same short exact sequence as in the proof of 
Lemma \ref{lm:ruleout2-2-1} that 
\begin{equation}   \label{eql:2-2-2-2dimcok-8}
h^0 (\mathrm{Cok}\, \gamma_1 \otimes \mathcal{O}_B(-8)) = 3.
\end{equation}

Since  
$\alpha_0 \beta_2 + \alpha_2 \beta_0 = 
\alpha_0 (c \alpha_0^2 - 2\beta_0)$, 
if the two sections $\alpha_0$ and $c \alpha_0^2 - 2\beta_0$ 
has a common zero $P$, then at this point $P$, the rank of 
$\sigma_2 \otimes k(P)$ drops at least by $2$, which is 
impossible. 
Thus $\mathrm{Cok}\, \gamma_1^{(1)}$ is locally free of 
rank $2$, where we denote by 
\begin{equation}   \label{eql:gamma_1^(1)}
\gamma_1^{(1)} : \mathcal{O}_B (5) \to 
\mathcal{O}_B (6)^{\oplus 2} \oplus \mathcal{O}_B (7)
\end{equation}
the multiplication morphism by 
$\,^t (-\alpha_0 , \, c\alpha_0, \, c \alpha_0^2 - 2\beta_0)$. 
Since the multiplication morphism 
$\,^t (-\alpha_0 ,\, c \alpha_0) \times : 
\mathcal{O}_B (5)  \to \mathcal{O}_B (6)^{\oplus 2}$ 
is the composite of the two morphisms $\alpha_0 \times : 
\mathcal{O}_B (5)  \to \mathcal{O}_B (6)$ 
and 
$\,^t (-1,\, c) \times : 
\mathcal{O}_B (6) \to \mathcal{O}_B (6)^{\oplus 2}$, 
we obtain by (\ref{eql:gamma_1^(1)}) that   
$h^0 (\mathrm{Cok}\, \gamma_1^{(1)} \otimes \mathcal{O}_B(-6)) = 4$,  
$h^0 (\mathrm{Cok}\, \gamma_1^{(1)} \otimes \mathcal{O}_B(-7)) = 2$, 
and 
$h^0 (\mathrm{Cok}\, \gamma_1^{(1)} \otimes \mathcal{O}_B(-8)) = 1$. 
From these together with 
$\deg \mathrm{Cok}\, \gamma_1^{(1)} = 14$, we infer  
$\mathrm{Cok}\, \gamma_1^{(1)} \simeq 
\mathcal{O}_B (6) \oplus \mathcal{O}_B (8)$. 
Thus by the same short exact sequence as in the proof of 
Lemma \ref{lm:ruleout2-2-1}, we see that 
$h^0 (\mathrm{Cok}\, \gamma_1 \otimes \mathcal{O}_B(-9)) = 1$, 
which together with (\ref{eql:2-2-2dimcok-6}) and 
(\ref{eql:2-2-2-2dimcok-8}) implies that 
$b_1 = 6$ and $b_i \geq 8$ for all $2 \leq i \leq 3$.  
Since $\deg \mathrm{Cok}\, \gamma_1 = 23$, we obtain 
\[
\mathrm{Cok}\, \gamma \simeq 
\mathcal{O}_B(4) \oplus \mathcal{O}_B(5) \oplus 
\mathcal{O}_B(6) \oplus
\mathcal{O}_B(8) \oplus \mathcal{O}_B(9),  
\]  
which contradicts Lemma \ref{lm:homl4cok}. \qed 

By Lemmas \ref{lm:homl4cok},  
\ref{lm:ruleout1}, 
\ref{lm:ruleout2-1-1}, 
\ref{lm:ruleout2-1-2}, 
\ref{lm:ruleout2-2-1}, 
\ref{lm:ruleout2-2-2-1}, and  
\ref{lm:ruleout2-2-2-2}, 
we obtain the following: 

\begin{proposition}  \label{prop:nevercompwap}
The map $\varPhi_{|2L|}$ is not composite with a pencil.    
\end{proposition}

{\bf Digression}. 
As we have already cited the result in our proof, 
Ciliberto--Francia--Mendes Lopes 
\cite[Proposition 1.7]{onbicanonicalmaps} shows 
that if a minimal surface $S$ of 
general type with $4 \leq p_g $ and $K^2 \leq 9$ has 
canonical map composite with a pencil $\mathcal{P}$, 
then either (i) $\mathcal{P}$ is a pencil of curves of 
genus $2$, or (ii) $\mathcal{P}$ is a rational pencil of 
non-hyperelliptic 
curves of genus $3$. 
Their proposition moreover shows that in the latter case 
$K^2 = 9$, $p_g= 4$, and $K \sim 3C$ hold, where 
$C$ is a general member of the pencil $\mathcal{P}$. 
As a byproduct of our computation, we prove the following: 

\begin{proposition}   \label{prop:nevoccurcfml}   
Case (ii) in \cite[Proposition 1.7]{onbicanonicalmaps} never occurs. 
Thus if a minimal surface $S$ with $4 \leq p_g $ and $K^2 \leq 9$ 
has canonical map composite with a pencil $\mathcal{P}$, 
then $\mathcal{P}$ is a pencil of curves of genus $2$. 
\end{proposition}

Proof. 
Assume that the surface $S$ has numerical invariants as 
in the assertion, and the canonical map $\varPhi_{|K|}$ is 
composite with a pencil $\mathcal{P}$ of curves of genus $3$.   
Then we have $K^2 = 9$, $p_g= 4$, and $K \sim 3C$, where 
$C$ is a general member of the pencil $\mathcal{P}$, which is 
non-hyperelliptic. 
If $S$ has irregularity $q = 0$, then it contradicts our 
Proposition \ref{prop:nevercompwap}. 
Thus it suffices to rule out the case $q > 0$, where 
$q$ is the irregularity of our surface $S$. 
In what follows we assume $q>0$. 

Note that we have $q \leq 2$. Indeed, 
since the restriction map 
$H^0 (\mathcal{O}_S (C)) \to H^0 (\mathcal{O}_C (C))$ has 
rank at least $1$, the map 
$H^0 (\mathcal{O}_S (4C)) \to H^0 (\mathcal{O}_C (4C))$
also has rank at least $1$. This together with the short exact 
sequence 
\[
0 \to \mathcal{O}_S(3C) = \mathcal{O}_S (K) 
\to \mathcal{O}_S (4C) \to \mathcal{O}_S (4C) = \omega_C \to 0
\]
implies the inequality.   
Thus we obtain $\chi = \chi (\mathcal{O}_S) \geq 3$. 
Then by so-called Severi inequality proved by Pardini, 
we see that 
$S$ is not of Albanese general type. 
We denote by $\alpha : S \to B$ the Albanese fibration 
of our surface $S$. Naturally we have $g(B) = q$, where $g(B)$ 
is the genus of the base curve $B$.  

First, let us rule out the case $q=2$. 
Assume that we have $q = 2$. Then $B$ is a non-singular curve 
of genus $g(B) = 2$. Since $\mathcal{P}$ has a unique 
base point $x$, the restriction 
$\alpha |_C : C \to B$ is surjective for any general 
$C \in \mathcal{P}$.   
Moreover, by Hurwitz formula, we see that $\alpha |_C$ is an 
\'{e}tale double cover. 
Let $F_0$ denote the fiber of $\alpha$ passing through 
the base point $x$. 
Then since $(F_0 \cdot C)_x =1$, we see that 
$x$ is a smooth point of $F_0$ and that 
for any $y \in F_0 \setminus \{ x \}$ 
there exists a unique member 
$C_y \in \mathcal{P}$ passing through $y$. 
This however contradicts the rationality of our 
pencil $\mathcal{P}$, since $F_0$ is a non-singular 
curve of genus $4$. 

Next, let us rule out the case $q=1$. 
We use the method used in the proof 
of \cite[Proposition 2.4]{onbicanonicalmaps}. 
For the reader's convenience, we include the outline of our proof. 
Assume that we have $q=1$. Then $B$ is an elliptic curve. 
We take a point $o \in B$, and use this point for the zero of 
the additive structure of the elliptic curve $B$. 
For any closed point $b \in B$, we put $\xi_b = \mathcal{O}_B (b -o)$. 

We first claim $h^0 (\mathcal{O}_C(C) \otimes \alpha^* \xi_b) =1$ 
for any $b \neq o \in B$.  
Indeed, assume otherwise. Then by the same method as in 
\cite[Proposition 2.4]{onbicanonicalmaps} 
and the upper-semicontinuity, we see that 
$h^0 (\mathcal{O}_S(C) \otimes \alpha^* \xi_{b^{\prime}}) = 0$ holds 
for any general $b^{\prime} \in B$. Then again by the same method 
as in \cite[Proposition 2.4]{onbicanonicalmaps}, we obtain 
$h^0 (\mathcal{O}_C(2C) \otimes \alpha^* \xi_{b^{\prime}}^{\vee}) \geq 3$ 
for any general $b^{\prime} \neq o \in C$. 
We however have $g(C) = 3$ and 
$\deg \mathcal{O}_C(2C) \otimes \alpha^* \xi_{b^{\prime}}^{\vee} = 2$. 
Thus this is impossible, since the curve  
$\varPhi_{|\mathcal{O}_C (2C) \otimes \alpha^* \xi_{b^{\prime}}^{\vee} |} (C)$ needs 
to be non-degenerate.  

Now put $U = B \setminus \{ o \}$, and denote by 
$\pi : C \times U \to U$ the second projection. 
Let $\varXi$ be an invertible sheaf on $C \times U$ such that 
$\varXi |_{C \times \{ b\}} \simeq 
\mathcal{O}_C (C) \otimes \alpha^* \xi_b$ holds for all $b \in U$.  
Then by what we have shown in the preceding paragraph, 
we see that the direct image $\pi_* \varXi$ is an invertible 
sheaf on $U$ 
and that the natural morphism 
$\pi^* \pi_* \varXi \to \varXi$ is non-trivial. 
Thus replacing $U$ by smaller one if necessary, 
we obtain an effective divisor 
$Z$ on $C \times U$ such that 
$Z \cap (C \times \{ b \}) = \mathrm{div}\, s_b$ holds  
for all $b \in U$, where 
$s_b \neq 0 \in H^0 (\mathcal{O}_C (C) \otimes \alpha^* \xi_b)$ 
is the unique non-zero global section to 
$\mathcal{O}_C (C) \otimes \alpha^* \xi_b$.  
Then the restriction $\pi |_Z : Z \to U$ is birational,  
since $\deg \mathcal{O}_C (C) \otimes \alpha^* \xi_b = 1$. 
By \cite[Proposition 1.6]{onbicanonicalmaps}, however, the first projection 
$Z \to C$ is dominant. Thus this contradicts the 
inequality $g(C)  =3 > g(B) =1$. \qed  

\section{Structure theorem}

Let us go back to the study of our surface $S$ with 
$c_1^2 = 9$, $\chi = 5$, and $K \sim 3L$.  
By Proposition \ref{prop:nevercompwap} 
and Lemma \ref{lm:notcompwap}, 
we have $h^0 (\mathcal{O}_S (L)) = 1$, 
$h^0 (\mathcal{O}_S (2L)) = 3$, 
and 
$2 \leq \deg \varPhi_{|2L|} \leq 4$, where $\deg \varPhi_{|2L|}$ 
is the degree of the rational map 
$\varPhi_{|2L|} : S - - \to \mathbb{P}^2$. 
In this section, we shall give a structure theorem for our 
surface $S$, by studying the structure of the graded ring 
$\bigoplus_{i=0}^{\infty} H^0(\mathcal{O}_S (nL))$.  

\begin{lemma}    \label{lm:fxdptcanonical}
Let $|K| = |3L| = |M_3| + F_3$ be the decomposition 
of the canonical system $|K|$ into the variable part 
$|M_3|$ and the fixed part $F_3$. 
Then $KF_3 = 0$ holds. In particular $F_3$ is at most 
a sum of fundamental cycles of rational double points. 
\end{lemma}

Proof. Let $M_3$ and $F_3$ be divisors as above. 
Let $p_3 : \tilde{S}_3 \to S$ be the shortest composite 
of quadric transformations such that the variable part 
$|\tilde{M}_3|$ of $p_3^* |M_3|$ is free from base points. 
Then we have $M_3^2 \geq \tilde{M}_3^2 \geq 4$ and 
$K^2 = M_3^2 + M_3 F_3 + K F_3$, where 
$M_3 F_3 \geq 0$ and $K F_3 \geq 0$ hold. 
Since we have $M_3^2 + M_3 F_3 = K M_3 \equiv 0 \mod 3$, 
this implies $M_3^2 + M_3 F_3 = 6$ or $9$. 

Assume that $M_3^2 + M_3 F_3 = 6$. 
Then by $M_3^2 \geq 4$, we have $0 \leq M_3 F_3 \leq 2$. 
From this together with Hodge's Index Theorem  
$M_3^2 F_3^2 = M_3^2 (K F_3 - M_3 F_3) \leq (M_3 F_3)^2$,  
we see that $M_3 F_3 = 2$, $M_3^2 =4$, and $F_3^2 =1$, 
hence $M_3 \sim_{\mathrm{num}} 2 F_3$ and 
$K \sim_{\mathrm{num}} 3 L \sim_{\mathrm{num}} 3 F_3$. 
Then by \cite[Theorem 4]{remarks'''''}, we obtain 
$L \sim F_3$ and $M_3 \sim 2 L$,    
which contradicts 
$h^0 (\mathcal{O}_S (2L)) = 3$ in Lemma \ref{lm:notcompwap},   
since we have $h^0 (\mathcal{O}_S (M_3)) = p_g(S) = 4$.  
Thus we obtain $M_3^2 + M_3 F_3 = 9$, hence the assertion. \qed   

Take a base $x_0$ of the space of global sections 
$H^0 (\mathcal{O}_S (L))$. 
The following lemma is trivial. 

\begin{lemma}  \label{lm:2L3LLstr}
$1$$)$ There exist two elements $y_0$, $y_1 \in H^0 (\mathcal{O}_S (2L))$   
such that $x_0^2$, $y_0$, and $y_1$ form a base of $H^0 (\mathcal{O}_S (2L))$. 

$2$$)$ There exists an element $z_0 \in H^0 (\mathcal{O}_S (3L))$ 
such that $x_0^3$, $x_0 y_0$, $x_0 y_1$, and $z_0$ form a base of 
$H^0 (\mathcal{O}_S (3L))$. 
\end{lemma}
  
Take three elements $x_0$, $y_1$, and $y_2$ as in the lemma above. 
In what follows. we denote by $C$ the unique member of the linear 
system $|L|$, and by $C_0$, its unique irreducible component 
such that $LC_0 = 1$. 
For the proof of the following lemma, see \cite[Lemma 1.2]{globalpg1k2=1}:  

\begin{lemma}   \label{lm:Dequalto2C}
If a member $D \in |2L|$ satisfies $D \succeq C_0$, then 
$D= 2C$. 
\end{lemma}

Let us study higher homogeneous parts of the ring 
$\bigoplus_{n=0}^{\infty} H^0 (\mathcal{O}_S (nL))$. 

\begin{lemma}     \label{lm:4Lstr}
The space $H^0 (\mathcal{O}_S (4L))$ has the decomposition 
$H^0 (\mathcal{O}_S (4L)) = 
x_0 H^0 (\mathcal{O}_S (3L)) \oplus 
\bigoplus_{i=0}^2 \mathbb{C}\, y_0^i y_1^{2-i}$. 
\end{lemma}

Proof. 
By the Riemann-Roch theorem, we have $h^0 (\mathcal{O}_S (4L)) = 7$. 
Thus it suffices to prove that seven elements 
$x_0^4$, $x_0^2 y_0$, $x_0^2 y_1$, $x_0 z_0$, 
$y_0^2$, $y_0 y_1$, and $y_1^2$ are linearly independent over 
$\mathbb{C}$. Assume that these seven elements has a nontrivial 
linear relation. Then there exist $(\alpha_0,\, \alpha_1)$ and 
$(\beta_0,\, \beta_1) \in \mathbb{C}^2 \setminus \{ 0 \} $ 
such that 
$(\alpha_0 y_0 + \alpha_1 y_1) (\beta_0 y_0 + \beta_1 y_1) 
\in x_0 H^0 (\mathcal{O}_S (3L))$.  
This contradicts Lemma \ref{lm:Dequalto2C}, since 
$x_0^2$, $y_0$, and $y_1$ are linearly independent. \qed

\begin{lemma}    \label{lm:5Lstr}   
There exists an element $u_0 \in H^0 (\mathcal{O}_S (5L))$ 
such that the equality 
$H^0 (\mathcal{O}_S (5L)) = 
x_0 H^0 (\mathcal{O}_S (4L))  
\oplus \bigoplus_{i=0}^1 \mathbb{C}\, y_0^i y_1^{1-i} z_0
\oplus \mathbb{C}\, u_0$ holds.   
\end{lemma}

Proof. By the Riemann-Roch theorem, we have 
$h^0 (\mathcal{O}_S (5L)) =10$. 
Thus it suffices to prove that a base of 
$x_0 H^0 (\mathcal{O}_S (4L))$ together with  
$y_0 z_0$ and $y_1 z_0$ forms a set of linearly 
independent nine elements of $H^0 (\mathcal{O}_S (5L)) $.  
Assume that these nine elements are not linearly independent 
over $\mathbb{C}$. Then there exists an element 
$(\alpha_0,\, \alpha_1) \in \mathbb{C}^2 \setminus \{ 0 \}$ 
such that 
$(\alpha_0 y_0 + \alpha_1 y_1) z_0 \in x_0 H^0 (\mathcal{O}_S (4L))$. 
The same argument as in the proof of Lemma \ref{lm:4Lstr} however     
shows that $\mathrm{div}\, (\alpha_0 y_0 + \alpha_1 y_1) \nsucceq C_0$. 
Thus we obtain $\mathrm{div}\, z_0 \succeq C_0$, which contradicts 
Lemma \ref{lm:fxdptcanonical}. \qed   

Take an element $u_0$ as in the lemma above. 
In what follows, we denote by 
$\mathbb{C} \left[ X_0, Y_0, Y_1, Z_0, U_0 \right]$ 
the weighted polynomial ring with 
$\deg X_0 = 1$, $\deg Y_0 = \deg Y_1 = 2$, 
$\deg Z_0 = 3$, and $\deg U_0 = 5$.  

\begin{lemma}     \label{lm:6Lstr}
There exists a homogeneous element  
$f_6 \in \mathbb{C} \left[ X_0, Y_0, Y_1, Z_0, U_0 \right]$ of 
degree $6$,  
unique up to multiplication by a non-zero constant, 
such that the equality 
$f_6 (x_0, y_0, y_1, z_0, u_0)  
= 0 \in H^0 (\mathcal{O}_S (6L))$ holds.  
The coefficient of $Z_0^2$ in $f_6$ is non-vanishing. 
Therefore the space of global sections to $6L$ decomposes as 
$H^0 (\mathcal{O}_S (6L)) = 
x_0 H^0 (\mathcal{O}_S (5L))  
\oplus \bigoplus_{i=0}^3 \mathbb{C}\, y_0^i y_1^{3-i}$. 
Moreover, by a proper choice of $z_0$, 
the polynomial $f_6$ can be set in such a way that it includes no 
term linear with respect to $Z_0$.  
\end{lemma}

Proof. 
The space $H^0 (\mathcal{O}_S (6L))$ contains $15$ monomials 
of $x_0$, $y_0$, $y_1$, $z_0$, and $u_0$; 
ten belonging to $x_0 H^0 (\mathcal{O}_S (5L))$, 
four of the form $y_0^i y_1^{3 - i}$ ($0 \leq i \leq 3$), 
and $z_0^2$. 
Meanwhile we have $h^0 (\mathcal{O}_S (6L)) = 14$. 
Thus there exists at least one non-trivial linear relation 
$f_6 (x_0, y_0, y_1, z_0, u_0) = 0$ among these $15$ monomials. 
Assume that the coefficient of $Z_0^2$ in $f_6$ vanishes. 
Then by Lemma \ref{lm:5Lstr}, there exist three elements 
$(\alpha_0,\, \alpha_1)$, $(\beta_0,\, \beta_1)$, 
$(\gamma_0,\, \gamma_1) \in \mathbb{C}^2 \setminus \{ 0 \}$ 
such that 
$(\alpha_0 y_0 + \alpha_1 y_1)
(\beta_0 y_0 + \beta_1 y_1)
(\gamma_0 y_0 + \gamma_1 y_1) 
\in x_0 H^0 (\mathcal{O}_S (5L))$, from which 
we infer a contradiction by the same argument 
as in Lemma \ref{lm:5Lstr}.     
Thus we obtain the non-vanishing of the 
coefficient of $Z_0^2$, hence also the 
uniqueness of $f_6$. 
The irreducibility of $f_6$ follows from 
Lemmas \ref{lm:2L3LLstr}, \ref{lm:4Lstr}, and \ref{lm:5Lstr}.   
\qed

In what follows, in view of the lemma above, we assume 
that $f_6$ includes no term linear with respect to $Z_0$. 

\begin{lemma}    \label{lm:bpfree2L}
The linear system $|2L|$ has no base point. 
Thus the map $\varPhi_{|2L|} : S  \to \mathbb{P}^2$ 
associated to $|2L|$ is a morphism of degree $4$. 
\end{lemma} 

Proof. 
Assume that the linear system $|2L|$ has a base point 
$b \in S$. Then this point $b$ is a common zero 
of $x_0$, $y_0$, and $y_1$. This together with Lemma \ref{lm:6Lstr} 
however implies that $b$ is a base point of $|6L| = |2K|$, 
which contradicts the base point freeness 
of the bicanonical system (see \cite[Theorem 2]{canonicalmdl}). 
\qed 
  
Now let  
$
\varphi_S : S \to 
\mathbb{P} (1, 2, 2, 3, 5) = 
\mathrm{Proj}\, \mathbb{C} \left[ X_0, Y_0, Y_1, Z_0, U_0 \right] 
$
denote the morphism induced by 
$X_0 \mapsto x_0$, $Y_i \mapsto y_i$ ($i = 0$, $1$), 
$Z_0 \mapsto z_0$, and $U_0 \mapsto u_0$. 

\begin{lemma}      \label{lm:varphiSbiratl}
The morphism $\varphi_S : S \to \mathbb{P} (1, 2, 2, 3, 5) $ 
is birational onto its image. 
\end{lemma}

Proof. Since the morphism $\varphi_S$ factors 
through the bicanonical map $\varPhi_{|2K|}$, 
it suffices to prove the birationality 
of $\varPhi_{|2K|} : S \to \mathbb{P}^{13}$. 
Assume that $\varPhi_{|2K|}$ is non-birational. 
Then by \cite[Theorem 1.8, Theorem 2.1]{onbicanonicalmaps}, 
the surface $S$ has a pencil $\mathcal{P}$ of curves of genus $2$. 
Moreover, by their proof, we see that $\mathcal{P}$ can be chosen 
in such a way that a general member $D \in \mathcal{P}$ satisfies 
$D^2 = 0$ and $DK= 2$. 
This however contradicts the equivalence $K \sim 3L$. \qed     

\begin{lemma}      \label{lm:7Lstr}
The space $H^0 (\mathcal{O}_S (7L))$ has the decomposition 
$H^0 (\mathcal{O}_S (7L)) = 
x_0 H^0 (\mathcal{O}_S (6L))  
\oplus \bigoplus_{i=0}^1 \mathbb{C}\, y_0^i y_1^{1-i} u_0 
\oplus 
\bigoplus_{i=0}^2 \mathbb{C}\, y_0^i y_1^{2-i} z_0$. 
\end{lemma}

Proof. 
Assume otherwise. 
Then, since $h^0 (\mathcal{O}_S (7L)) = 19$ and 
$h^0 (\mathcal{O}_S (6L)) $$= 14$, 
there exists a non-zero homogeneous element 
$g_7 \in \mathbb{C} \left[ X_0, Y_0, Y_1, Z_0, U_0 \right]$ 
of degree $7$  
satisfying $g_7 (x_0, y_0, y_1, z_0, u_0) = 0 \in H^0 (\mathcal{O}_S (7L))$ 
in which at least one of the five monomials   
$Y_0 U_0$, $Y_1 U_0$, $Y_0^2 Z_0$, $Y_0 Y_1 Z_0$, and $Y_1^2 Z_0$ 
has non-vanishing coefficient.  
By subtracting a multiple of $f_6$, we may assume that 
the coefficient of $X_0 Z_0^2$ in $g_7$ vanishes. 
Moreover, $Y_0 U_0$ or $Y_1 U_0$ has 
non-vanishing coefficient in $g_7$. 
Indeed, if both $Y_0 U_0$ and $Y_1 U_0$ have vanishing coefficient, 
then the same argument as in Lemma \ref{lm:5Lstr} shows that 
$\mathrm{div}\, z_0 \succeq C_0$, which contradicts 
Lemma \ref{lm:fxdptcanonical}.    
 
Now let $\mathcal{Q}_7 \subset \mathbb{P} (1, 2, 2, 3, 5)$ 
be the subvariety defined by $f_6 = g_7 =0$, 
and $\pi_{\mathcal{Q}_7} : \mathcal{Q}_7 - - \to \mathbb{P} (1, 2, 2)$, 
the restriction to $\mathcal{Q}_7$ of the natural dominant map 
$\mathbb{P} (1, 2, 2, 3, 5) - - \to \mathbb{P} (1, 2, 2)$. 
Since $g_7$ is not a multiple of $f_6$ in 
$\mathbb{C} \left[ X_0, Y_0, Y_1, Z_0, U_0 \right]$, 
and since $f_6$ is irreducible, we have $\dim \mathcal{Q}_7 = 2$. 
Moreover, since $Z_0$ appears quadratically and $U_0$ appears 
at most linearly in $f_6$, 
and since $U_0$ appears linearly and $Z_0$ appears at most 
linearly in $g_7$, we see that 
$\deg \pi_{\mathcal{Q}_7} \leq 2$.  
Meanwhile, since 
$\pi_{\mathcal{Q}_7} \circ \varphi_S : S \to \mathbb{P} (1 ,2, 2)$ 
coincides with $\varPhi_{|2L|}$ via the natural isomorphism 
$\mathbb{P} (1 ,2, 2) \simeq 
\mathbb{P}^2 = \mathbb{P} (H^0 (\mathcal{O}_S (2L)))$, 
we see by Lemma \ref{lm:bpfree2L} that 
$\deg \pi_{\mathcal{Q}_7} \circ \varphi_S = 4$. 
Thus we obtain $\deg \varphi_S \geq 2$, which contradicts 
Lemma \ref{lm:varphiSbiratl}.  \qed 

By the same method, we can prove the following two lemmas:

\begin{lemma}
The space $H^0 (\mathcal{O}_S (8L))$ has the decomposition 
$H^0 (\mathcal{O}_S (8L)) = 
x_0 H^0 (\mathcal{O}_S (7L))  
\oplus \mathbb{C}\, z_0 u_0 
\oplus 
\bigoplus_{i=0}^4 \mathbb{C}\, y_0^i y_1^{4-i}$. 
\end{lemma}

\begin{lemma}      \label{lm:9Lstr}
The space $H^0 (\mathcal{O}_S (9L))$ has the decomposition 
$H^0 (\mathcal{O}_S (9L)) = 
x_0 H^0 (\mathcal{O}_S (8L))  
\oplus 
\bigoplus_{i=0}^2 \mathbb{C}\, y_0^i y_1^{2-i} u_0
\oplus 
\bigoplus_{i=0}^3 \mathbb{C}\, y_0^i y_1^{3-i} z_0$. 
\end{lemma}

Indeed, 
we just need to consider 
$\mathcal{Q}_k = \{ f_6 = g_k =0 \} \subset 
\mathbb{P} (1, 2, 2, 3, 5)$ and 
$\pi_{\mathcal{Q}_k} : \mathcal{Q}_k - - \to \mathbb{P} (1, 2, 2)$ 
for $k = 8$, $9$: we obtain easily   
$\deg \pi_{\mathcal{Q}_k} \leq 3$ for $k = 8$, $9$, 
which leads us to a contradiction to 
Lemma \ref{lm:varphiSbiratl}.

\begin{corollary}      \label{cor:5Lbpfree}
The linear system $|5L|$ has no base point. 
\end{corollary}

Proof. 
Note that by Lemma \ref{lm:bpfree2L} the three sections 
$x_0^2$, $y_0$, $y_1 \in H^0 (\mathcal{O}_S (2L))$ have 
no common zero. 
Since we have 
$x_0^5$, $x_0^2z_0$, $y_0z_0$, $y_1z_0$, 
$u_0 \in H^0 (\mathcal{O}_S (5L))$, 
this implies that the base locus of $|5L|$ 
is contained in the subset  
$\{ x_0 = z_0 = u_0 = 0 \} \subset S$.  
Thus, by Lemma \ref{lm:9Lstr}, 
if the linear system $|5L|$ has a base point $b \in S$, 
then this point $b$ is also a base point of 
$|9L| = |3K|$, which contradicts the  
base point freeness of the tricanonical system 
(see \cite[Theorem 2]{canonicalmdl}). \qed

\begin{lemma}
There exists a homogeneous element 
$g_{10} \in \mathbb{C} \left[ X_0, Y_0, Y_1, Z_0, U_0 \right]$ 
of degree $10$ not multiple of $f_6$ 
such that 
$g_{10} (x_0, y_0, y_1, z_0, u_0) =0$ holds in 
$H^0 (\mathcal{O}_S (10L))$.  
The coefficient of $U_0^2$ in $g_{10}$ is non-vanishing. 
The polynomial $g_{10}$ can be chosen in such a way that 
it includes no monomial divisible by $Z_0^2$, and 
with this last condition imposed, the polynomial    
$g_{10}$ is unique up to 
multiplication by a non-zero constant. 
Moreover the space  
$H^0 (\mathcal{O}_S (10L))$ 
decomposes as 
$H^0 (\mathcal{O}_S (10L)) = 
x_0 H^0 (\mathcal{O}_S (9L))  
\oplus \bigoplus_{i=0}^1 \mathbb{C}\, y_0^i y_1^{1-i} z_0 u_0 
\oplus 
\bigoplus_{i=0}^5 \mathbb{C}\, y_0^i y_1^{5-i}$.     
\end{lemma}

Proof. 
The space $H^0 (\mathcal{O}_S (10L))$ includes $41$ monomials 
not divisible by $z_0^2$ of 
$x_0$, $y_0$, $y_1$, $z_0$, and $u_0$. 
Since $h^0 (\mathcal{O}_S (10L)) = 40$, this implies that 
there exists at least one homogeneous element 
$g_{10} \in \mathbb{C} \left[ X_0, Y_0, Y_1, Z_0, U_0 \right]_{10}$ 
as in the first assertion. 
Since $g_{10}$ includes no monomial divisible by $Z_0^2$, 
it is not a multiple of $f_6$. 
Assume that the coefficient of $U_0^2$ in $g_{10}$ vanishes. 
Then by the same argument as in the proof 
of Lemma \ref{lm:5Lstr}, we see that 
ether the coefficient of $Y_0 Z_0 U_0$ or that of 
$Y_1 Z_0 U_0$ is non-vanishing. 
This however together with the same argument 
as in the proof of 
Lemma \ref{lm:7Lstr} leads us to a contradiction 
to Lemma \ref{lm:varphiSbiratl}.   
Thus the coefficient of $U_0^2$ is non-vanishing, from which  
the last assertion and the uniqueness of $g_{10}$ follow.  \qed 

Let $\mathcal{Q}$ denote the subvariety of 
$\mathbb{P} (1, 2, 2, 3, 5)$ defined by the ideal  
$(f_6 , g_{10})$. 
We define the subvarieties 
$\mathcal{Z}_0$, $\mathcal{Z}_1$, and $\mathcal{Z}_2$ 
of $\mathbb{P} (1, 2, 2, 3, 5)$ by   
\begin{align}
\mathcal{Z}_0 &= \{ X_0 = Z_0 = U_0 =0 \}, \notag \\ 
\mathcal{Z}_1 &= \{ X_0 = Y_0 = Y_1 = U_0 =0 \},  \notag \\  
\mathcal{Z}_2 &= \{ X_0 = Y_0 = Y_1 = Z_0 =0 \}.& &  \notag
\end{align}
Note that outside $\bigcup_{i=0}^2 \mathcal{Z}_i$ the 
weighted projective space $\mathbb{P} (1, 2, 2, 3, 5)$ 
has no singularity. 
The restriction of $\mathcal{O} (1)$ to 
$\mathbb{P} (1, 2, 2, 3, 5) \setminus \bigcup_{i=0}^2 \mathcal{Z}_i$ 
is invertible. 

\begin{proposition}
$1$$)$  
The morphism $\varphi_S : S \to \mathbb{P} (1, 2, 2, 3, 5)$
surjects to $\mathcal{Q}$.  

$2$$)$ The variety $\mathcal{Q}$ does not intersect the locus  
$\bigcup_{i=0}^2 \mathcal{Z}_i$. 

$3$$)$ The inclusion map
$\varphi_S^* : \mathbb{C} \left[ X_0, Y_0, Y_1, Z_0 , U_0 \right] / 
(f_6, g_{10}) \to R(S, L)$ 
is an isomorphism of graded $\mathbb{C}$-algebra, where 
$R(S, L): = \bigoplus_{n=0}^{\infty} H^0 (\mathcal{O}_S (nL))$. 
The variety $\mathcal{Q}$ has at most rational double points 
as its singularities. 
\end{proposition}

Proof. 
Since $\deg \mathcal{Q} = 1 = L^2$, the assertion $1$) follows 
from Lemma \ref{lm:varphiSbiratl}. 
Then the assertion $2$) follows from the non-vanishing of 
the coefficient of $Z_0^2$ in $f_6$, 
that of the the coefficient of $U_0^2$ in $g_{10}$, 
Lemma \ref{lm:5Lstr}, and Corollary \ref{cor:5Lbpfree}. 
It only remains to prove the assertion $3$). 
By the assertion $2$), we see that $\mathcal{Q}$ is 
Gorenstein. Moreover we have 
$\omega_{\mathcal{Q}} \simeq \mathcal{O}_{\mathcal{Q}} (3)$, 
hence $\omega_S \simeq \varphi_S^* \omega_{\mathcal{Q}}$. 
Thus $\mathcal{Q}$ has at most rational double points 
as its singularities. 
Since $\varphi_S : S \to \mathcal{Q}$ 
gives the minimal desingularization 
of $\mathcal{Q}$, we obtain the assertion $3$). \qed    

Naturally,  
$R (S, L)^{(3)} = \bigoplus_{n=0}^{\infty} H^0 (\mathcal{O}_S (3nL))$ 
is the canonical ring of the surface $S$. 
Thus, we obtain the following: 

\begin{theorem}      \label{thm:structurethm}
If a minimal surface $S$ has $c_1^2 = 9$ and $\chi =5$, and  
its canonical class is divisible by $3$ in its integral cohomology 
group,  
then its canonical model is a $(6, 10)$-complete intersection 
of the weighted projective space $\mathbb{P} (1, 2, 2, 3 ,5)$ 
that does not intersect the locus $\bigcup_{i=0}^2 \mathcal{Z}_i$. 
Conversely, if a $(6, 10)$-complete intersection 
$\mathcal{Q} \subset \mathbb{P} (1, 2, 2, 3 ,5)$ satisfying     
$\mathcal{Q} \cap \bigcup_{i=0}^2 \mathcal{Z}_i = \emptyset$ 
has at most rational double points as its singularities, 
then its minimal desingularization $S$ is a minimal surface 
with $c_1^2 = 9$ and $\chi =5$ whose canonical 
class is divisible by $3$.           
\end{theorem} 

 
Note that for general $f_6$ and 
$g_{10} \in \mathbb{C} \left[ X_0, Y_0, Y_1, Z_0, U_0 \right]$ 
of degree $6$ and $10$, respectively, 
the subvariety $\mathcal{Q} = \{ f_6 = g_{10} = 0 \}
\subset \mathbb{P} (1, 2, 2, 3, 5)$ is non-singular. 
This can be verified with 
$X_0^6$, $Y_0^3$, $Y_1^3$, 
$Z_0^2 \in \mathbb{C} \left[ X_0, Y_0, Y_1, Z_0, U_0 \right]_6$,    
$X_0^{10}$, $Y_0^5$, $Y_1^5$, 
$U_0^2 \in \mathbb{C} \left[ X_0, Y_0, Y_1, Z_0, U_0 \right]_{10}$,   
and Bertini's Theorem.

\begin{remark}
Let $S$ and $S^{\prime}$ be two minimal algebraic surfaces 
with invariants as in Theorem \ref{thm:structurethm}. 
Then as one can see from the proof of Theorem \ref{thm:structurethm}, 
the surfaces  
$S$ and $S^{\prime}$ are isomorphic to each other, if and only if  
the varieties $\mathcal{Q}$ and $\mathcal{Q}^{\prime}$ are projectively 
equivalent 
in the weighted projective space 
$\mathbb{P} (1, 2, 2, 3, 5)$, where 
$\mathcal{Q}$ and $\mathcal{Q}^{\prime}$ are 
the $(6, 10)$-complete intersections corresponding to 
$S$ and $S^{\prime}$, respectively.    
\end{remark}

\begin{remark}       \label{rmk:eqcdforeptness}
Let $f_6$ and $g_{10} \in \mathbb{C} \left[ X_0, Y_0, Y_1, Z_0, U_0 \right]$ 
be homogeneous polynomials of weighted degree $6$ and $10$, respectively.  
Assume that the coefficient of $Z_0^2$ in $f_6$ 
and that of $U_0^2$ in $g_{10}$ are non-vanishing. 
Let $\mathcal{Q} \subset \mathbb{P} (1, 2, 2, 3, 5)$ denote the 
subvariety defined by the polynomials $f_6$ and $g_{10}$. 
Then  
$\mathcal{Q} \cap \bigcup_{i=0}^2 \mathcal{Z}_i = \emptyset$ 
holds, if and only if the two sections 
$f_6 (0, Y_0, Y_1, 0, 0) \in H^0 (\mathcal{O}_{\mathbb{P}^1} (3))$ 
and 
$g_{10} (0, Y_0, Y_1, 0, 0) \in H^0 (\mathcal{O}_{\mathbb{P}^1} (5))$ 
have no common zero on the projective line 
$\mathbb{P}^1 = \mathrm{Prpj}\, \mathbb{C} \left[ Y_0, Y_1 \right]$.      
\end{remark}

\section{Moduli space and the canonical maps} 

In this section, we study the moduli space. 
We also study the behavior of the canonical map of our surface $S$. 
Let us begin with the normal form of the defining polynomials 
$f_6$ and $g_{10}$.

\begin{proposition}    \label{prop:normalform}
Let $S$ be a minimal surface as in Theorem \ref{thm:structurethm}. 
Then the defining polynomials $f_6$ and $g_{10}$ in 
$\mathbb{P} (1, 2, 2, 3, 5)$ of its canonical model $\mathcal{Q}$ 
can be taken in the form
\begin{align}
f_6 &= Z_0^2 + \alpha_0 X_0 U_0 + \alpha_3 (X_0^2, Y_0, Y_1),& \notag \\
g_{10} &= 
U_0^2 + \beta_3 (X_0^2, Y_0, Y_1) X_0 Z_0 + \beta_5 (X_0^2, Y_0, Y_1),& \notag
\end{align}
where $\alpha_0 \in \mathbb{C}$ is a constant, 
$\alpha_3$, a homogeneous polynomial of degree $3$,  
and $\beta_i$, a homogeneous polynomial of degree $i$ 
for $i = 3$, $5$.  
\end{proposition}

Proof. 
By completing the square with respect to $Z_0$, we can take 
$f_6$ and $g_{10}$ in the form  
\begin{align}
f_6 &= Z_0^2 + \alpha_0 X_0 U_0 + \alpha_3 (X_0^2, Y_0, Y_1),& \notag \\
g_{10} &= 
U_0^2 + \beta_1 (X_0^2, Y_0, Y_1) Z_0 U_0 + 
\beta_3 (X_0^2, Y_0, Y_1) X_0 Z_0 + \beta_5 (X_0^2, Y_0, Y_1).& \notag
\end{align}
Putting $X_0 = X_0^{\prime}$, $Y_0 = Y_0^{\prime}$, $Y_1 = Y_1^{\prime}$, 
$Z_0 = Z_0^{\prime} + \alpha_0 \beta_1 X_0^{\prime} / 4 $, and    
$U_0 = U_0^{\prime} - \beta_1 Z_0^{\prime} / 2  
- \alpha_0 \beta_1^2 X_0^{\prime} / 4 $, and 
employing $X_0^{\prime}$, $Y_0^{\prime}$, $Y_1^{\prime}$, 
$Z_0^{\prime}$, $U_0^{\prime}$ as new 
$X_0$, $Y_0$, $Y_1$, 
$Z_0$, $U_0$, respectively,  
we easily obtain new $f_6$ and $g_{10}$ in which  
the term $\beta_1 Z_0 U_0$ vanishes. \qed    

Using this proposition, we prove the following theorem:

\begin{theorem}     \label{thm:modulisp}
The coarse moduli space $\mathcal{M}$ of surfaces  as in 
Theorem \ref{thm:structurethm} is a unirational variety 
of dimension $34$. 
In particular, any two surfaces $S$'s as in 
Theorem \ref{thm:structurethm}   
are deformation equivalent to each other. 
\end{theorem}

Proof. 
In what follows, for two weighted homogeneous 
polynomials $f_6$ and $g_{10}$ as 
in Proposition \ref{prop:normalform}, 
we denote by $S_{(f_6 ,\, g_{10})}$ the minimal desingularization 
of the variety $ \mathcal{Q}_{(f_6 ,\, g_{10})} = 
\{ f_6 = g_{10} = 0 \} \subset \mathbb{P} (1, 2, 2, 3, 5)$.  
Note that the pair $(f_6 ,\, g_{10})$ in the normal form 
as in Proposition \ref{prop:normalform} has $42$ linear parameters.  
Denote by $V$ the Zariski open subset 
of $\mathbb{A}^{42}$ consisting of all      
$(f_6, \, g_{10})$'s such that 
1) $\mathcal{Q}_{(f_6 ,\, g_{10})}$ has at most rational double points 
as its singularities, and 
2) $\mathcal{Q}_{(f_6 ,\, g_{10})} \cap \bigcup_{i=0}^2 \mathcal{Z}_i
= \emptyset$ holds. 
Then by the existence of the natural family of the canonical models 
$\mathcal{Q}_{(f_6 ,\, g_{10})}$'s over the space of parameters $V$,  
we obtain the irreducibility of the moduli space $\mathcal{M}$. 

Let us compute the dimension of the moduli space $\mathcal{M}$. 
Note that for two points $(f_6,\, g_{10})$ and
 $(f_6^{\prime},\, g_{10}^{\prime})$  of $V$, 
the corresponding surfaces 
$S_{(f_6,\, g_{10})}$ and $S_{(f_6^{\prime},\, g_{10}^{\prime})}$ are isomorphic 
to each other if and only if the ideals $(f_6,\, g_{10})$ and  
$(f_6^{\prime},\, g_{10}^{\prime})$ 
are equivalent under the action by the 
group of homogeneous transformations on the graded algebra  
$\mathbb{C} \left[ X_0, Y_0, Y_1, Z_0, U_0 \right]$. 
Since no monomial divisible by $Z_0^2$ appears  
in $g_{10}$ and $g_{10}^{\prime}$, 
and since $10 - \deg U_0 < \deg Z_0^2$ holds, this is equivalent to 
the condition that the points $(f_6,\, g_{10})$ and  
$(f_6^{\prime},\, g_{10}^{\prime})$ of $V$ are equivalent under the action 
by the group of homogeneous transformations of  
$\mathbb{C} \left[ X_0, Y_0, Y_1, Z_0, U_0 \right]$.   
Moreover, we see easily that if a point $v\in V$ corresponding to a 
surface $S$ is sufficiently general, 
then there exists a point $(f_6,\, g_{10}) \in V$ 
that gives the same isomorphism class 
of $S$ and such that $f_6$ and $g_{10}$ are in the form
\begin{align}
f_6 &= Z_0^2 + X_0 U_0 + Y_0^3 + Y_1^3 + a_0 X_0^2 Y_0Y_1 
+ X_0^4 (a_1 Y_0 + a_2 Y_1) + X_0^6,& \notag \\
g_{10} &= 
U_0^2 +  
\beta_3 (X_0^2, Y_0, Y_1) X_0 Z_0 + \beta_5 (X_0^2, Y_0, Y_1).& \notag
\end{align}
We denote by $V^{\prime}$ the $34$-dimensional subvariety of $V$ 
consisting of all $(f_6, \, g_{10})$'s in this form. 
Then the restriction $V^{\prime} \to \mathcal{M}$ of 
the natural morphism $V \to \mathcal{M}$ is dominant.  

Let us study fibers of the morphism 
$V^{\prime} \to \mathcal{M}$.   
Let $G$ be the group of homogeneous transformations of 
$\mathbb{C} \left[ X_0, Y_0, Y_1, Z_0, U_0 \right]$ that 
preserve the subvariety $V^{\prime} \subset V$. 
We denote by $\omega$ the third root of unity, and define 
the two transformations $\sigma$, $\tau \in G$ by   
\begin{align}
\sigma &:& U_0 &\mapsto  U_0,& Z_0 &\mapsto Z_0& 
         Y_0 &\mapsto \omega Y_0,&  Y_1 &\mapsto \omega^2 Y_1& 
         X_0 &\mapsto  X_0  & \notag \\  
\tau   &:& U_0 &\mapsto  U_0,& Z_0 &\mapsto Z_0& 
         Y_0 &\mapsto Y_1,&  Y_1 &\mapsto Y_0& 
         X_0 &\mapsto  X_0 . & \notag  
\end{align}
Then $\sigma$ and $\tau$ generate a subgroup 
$\langle \sigma, \tau \rangle \simeq \mathfrak{S}_3 \subset G$, 
where $\mathfrak{S}_3$ is the symmetric group of degree $3$. 
Since each element of $G$ induces a permutation of three prime divisors of 
$Y_0^3 + Y_1^3$, we have a natural group homomorphism 
$G \to \mathfrak{S}_3$, whose restriction   
$\langle \sigma, \tau \rangle \to \mathfrak{S}_3$ to 
$\langle \sigma, \tau \rangle \subset G$ is an isomorphism. 
Thus if we define  
$\varPsi_{(\lambda_0 , \lambda_1 , \mu_0 , a )} \in G$ by 
$U_0 \mapsto a^5 U_0$,  
$Z_0 \mapsto (-1)^{\mu_0} a^3 Z_0$, 
$Y_1 \mapsto \omega^{\lambda_1} a^2 Y_1$, 
$Y_0 \mapsto \omega^{\lambda_0} a^2 Y_0$, and 
$X_0 \mapsto a X_0$ 
for each 
$(\lambda_0 , \lambda_1 , \mu_0 , a ) \in 
(\mathbb{Z} / 3)^{\oplus 2} \oplus 
\mathbb{Z} / 2 \oplus \mathbb{C}^{\times}$, 
then each element of $G$ can be written as 
$\rho \circ \varPsi_{(\lambda_0 , \lambda_1 , \mu_0 , a )}$ 
for an element $\rho \in \langle \sigma , \tau \rangle$ 
and an element 
$(\lambda_0 , \lambda_1 , \mu_0 , a ) \in 
(\mathbb{Z} / 3)^{\oplus 2} \oplus 
\mathbb{Z} / 2 \oplus \mathbb{C}^{\times}$. 
This implies that for a general point of $\mathcal{M}$   
the fiber of $V^{\prime} \to \mathcal{M}$ over this point consists 
of at most $108$ points. 
Now since $V^{\prime}$ is a Zariski open subset of the affine space  
$\mathbb{A}^{34}$, we see that $\mathcal{M}$ is unirational of 
dimension $34$. 
\qed    

Finally, we study the behavior of the canonical map of our surface $S$. 
Let $\mathcal{Q} \simeq 
\mathrm{Proj}\, 
\mathbb{C} \left[ X_0, Y_0, Y_1, Z_0, U_0 \right] / (f_6, g_{10})$  
be the canonical model of our surface $S$, where 
$f_6$ and $g_{10}$ are in the normal form as in 
Proposition \ref{prop:normalform}. 
Since the birational morphism $\varphi_S : S \to \mathcal{Q}$ 
factors through the canonical map $\varPhi_{|K|} : S - - \to \mathbb{P}^3$, 
the study of the behavior of $\varPhi_{|K|}$ is reduced to 
that of the behavior of the rational map 
$\varPhi_{|\mathcal{O}_{\mathcal{Q}} (3)|} : 
\mathcal{Q} - - \to \mathbb{P}^3$.  
Let $\xi_0$, $\eta_0$, $\eta_1$, and $\zeta_0$ be the 
homogeneous coordinates of $\mathbb{P}^3$ corresponding to 
the base $X_0^3$, $X_0 Y_0$, $X_0 Y_1$, $Z_0$ of 
$H^0 (\mathcal{O}_{\mathcal{Q}} (3))$. 
Note that for an integer $d \geq 1$, an equation of 
$\varPhi_{|K|} (S)$ in $\mathbb{P}^3$ of degree $d$ 
corresponds to a relation among 
$X_0^3$, $X_0 Y_0$, $X_0 Y_1$, and $Z_0$ in the 
homogeneous part of degree $3d$ of  
$\mathbb{C} \left[ X_0, Y_0, Y_1, Z_0, U_0 \right] / (f_6, g_{10})$. 

\begin{theorem}    \label{thm:canonicalimage}
Let $S$ be a minimal surface as in Theorem \ref{thm:structurethm}, 
and $f_6$ and $g_{10}$, the defining polynomials in 
$\mathbb{P} (1, 2, 2, 3, 5)$ of its canonical model $\mathcal{Q}$. 
Assume that $f_6$ and $g_{10}$ are in the normal form as in 
Proposition \ref{prop:normalform}.  

$1$$)$ If $\alpha_0 \neq 0$, then the canonical map 
$\varPhi_{|K|}$ of $S$ is birational onto its image, and 
the canonical image $\varPhi_{|K|} (S)$ is a sextic surface 
in $\mathbb{P}^3$ defined by 
\[
\left[ \xi_0 \zeta_0^2 + \alpha_3 (\xi_0 , \eta_0 , \eta_1) \right]^2 
+ \alpha_0^2 \left[ \beta_3 (\xi_0 , \eta_0 , \eta_1) \xi_0^2 \zeta_0 
+ \beta_5 (\xi_0 , \eta_0 , \eta_1) \xi_0 \right] = 0.
\]
Surfaces $S$'s with birational $\varPhi_{|K|}$ form an 
open dense subset of $\mathcal{M}$.  

$2$$)$ If $\alpha_0 = 0$, then the canonical map 
$\varPhi_{|K|}$ of $S$ is generically two-to-one onto its image, 
and the canonical image $\varPhi_{|K|} (S)$ is a cubic surface 
in $\mathbb{P}^3$ defined by 
\[
\xi_0 \zeta_0^2 + \alpha_3 (\xi_0 , \eta_0 , \eta_1) = 0.
\]
Surfaces $S$'s with non-birational $\varPhi_{|K|}$ form 
a $33$-dimensional locus in $\mathcal{M}$.  
\end{theorem}

Proof. 
The only non-trivial relation in the homogeneous part 
of degree $6$ of 
$\mathbb{C} \left[ X_0, Y_0, Y_1, Z_0, U_0 \right] / (f_6, g_{10})$  
is given by  $f_6 = Z_0^2 + \alpha_0 X_0 U_0 + \alpha_3 (X_0^2, Y_0, Y_1) = 0$. 
Assume that $f_6$ is a polynomial of 
$X_0^3$, $X_0 Y_0$, $X_0 Y_1$ and $Z_0$. 
Then $\alpha_3 (0, Y_0 , Y_1)$ must be zero in 
$\mathbb{C} \left[ Y_0, Y_1 \right]$. 
In this case, however, we have 
$f_6 (0, Y_0, Y_1, 0, 0) = 0 \in \mathbb{C} \left[ Y_0, Y_1 \right]$, 
which contradicts the condition 
$\mathcal{Q} \cap \bigcup_{i=0}^2 \mathcal{Z}_i = \emptyset$ 
(see Remark \ref{rmk:eqcdforeptness}). 
Thus $\varPhi_{|K|} (S) \subset \mathbb{P}^3$ satisfies no equation 
of degree $2$.  

Assume that $X_0^3$, $X_0 Y_0$, $X_0 Y_1$, and $Z_0$ have a 
non-trivial relation in the homogeneous part of degree $9$ of 
$\mathbb{C} \left[ X_0, Y_0, Y_1, Z_0, U_0 \right] / (f_6, g_{10})$. 
Then this relation must be written as 
$\gamma_1 (X_0^3, X_0 Y_0, X_0 Y_1, Z_0) f_6 = 0$, 
where $\gamma_1$ is a linear form with coefficients in $\mathbb{C}$. 
Since this left hand is a polynomial of 
$X_0^3$, $X_0 Y_0$, $X_0 Y_1$, and $Z_0$, 
we see with the help of Remark \ref{rmk:eqcdforeptness} 
that $\alpha_0 = 0$ holds and that 
$\gamma_1 (X_0^3, X_0 Y_0, X_0 Y_1, Z_0)$ is a multiple of $X_0^3$. 
Thus if $\alpha_0 = 0$, then $\varPhi_{|K|} (S) \subset \mathbb{P}^3$ is 
a cubic surface as in the assertion, and if $\alpha_0 \neq 0$, then 
$\varPhi_{|K|} (S) \subset \mathbb{P}^3$ satisfies no equation of 
degree $3$.  

Now that we have shown the assertion for the case $\alpha_0 = 0$,  
we assume in what follows that $\alpha_0 \neq 0$. 
By an argument similar to that in the preceding paragraph, 
we can prove the absence of equations 
of degree $d$ of $\varPhi_{|K|} (S) \subset \mathbb{P}^3$  
for $d = 4$, $5$. 
On the other hand, we can easily find an equation of degree $6$ that is 
satisfied by $\varPhi_{|K|} (S) \subset \mathbb{P}^3$.   
Note that in 
$\mathbb{C} \left[ X_0, Y_0, Y_1, Z_0, U_0 \right]$ we have 
$ -\alpha_0 X_0 U_0 \equiv Z_0^2 + \alpha_3 (X_0^2, Y_0, Y_1)$ 
and 
$ -U_0^2 \equiv 
\beta_3 (X_0^2, Y_0, Y_1) X_0 Z_0 + \beta_5 (X_0^2, Y_0, Y_1)$ 
modulo the ideal $(f_6, g_{10})$.   
Eliminating $U_0$ from these two and then multiplying it by $X_0^6$, 
we obtain 
\begin{multline}
\left[ X_0^3 Z_0^2 + \alpha_3(X_0^3, X_0 Y_0, X_0 Y_1) \right]^2   \notag  \\
+ \alpha_0 X_0^3 
\left[ \beta_3(X_0^3, X_0 Y_0, X_0 Y_1) X_0^3 Z_0 
 +\beta_5(X_0^3, X_0 Y_0, X_0 Y_1) \right]  \equiv 0 \notag 
\end{multline}
modulo the ideal $(f_6, g_{10})$. 
From this together with the absence of equation of lower degree, 
we see that if $\alpha_0 \neq 0$ then 
$\varPhi_{|K|} (S) \subset \mathbb{P}^3$ is a sextic surface 
defined by the equation as in the assertion. 

Now let us compute the mapping degree of the canonical map 
$\varPhi_{|K|} : S - - \to \mathbb{P}^3$.  
Let $|K| = |3L| = |M_3| + F_3$ be the decomposition 
as in Lemma \ref{lm:fxdptcanonical}, and 
$p_3 : \tilde{S}_3 \to S$, the shortest composite of 
quadric transformations such that the variable part 
$|\tilde{M}_3|$ of $p_3^* |M_3|$ is free from base points.     
Then we have 
\begin{equation}   \label{eql:mappingdeg}
\deg \varPhi_{|K|} \deg \varPhi_{|K|} (S) =
\tilde{M}_3^2 \leq M_3^2 \leq K^2 = 9. 
\end{equation}
Assume that $\alpha_0 \neq 0$. Then since 
$\deg \varPhi_{|K|} (S) =6$, we infer from the inequalities above  
that $\deg \varPhi_{|K|} = 1$.  
Assume that $\alpha_0 = 0$. Then since 
$\deg \varPhi_{|K|} (S) =3$, we infer in the same way 
that $\deg \varPhi_{|K|} \leq 3$. 
If $\deg \varPhi_{|K|} = 3$ holds, however, we see from (\ref{eql:mappingdeg}) 
that the linear system $|K|$ is base point free. 
This is impossible, because by Lemma \ref{lm:2L3LLstr} the canonical 
system $|K|$ needs to have a base point. 
Thus we obtain $\deg \varPhi_{|K|} = 2$.  

It is trivial that the surfaces $S$'s with birational 
$\varPhi_{|K|}$ form 
an open dense subset in $\mathcal{M}$. 
To show that the surfaces $S$'s with non-birattional 
$\varPhi_{|K|}$ form 
a $33$-dimensional locus in $\mathcal{M}$, we just need to 
use the same method as in the computation of $\deg \mathcal{M}$
in Theorem \ref{thm:modulisp}.  
\qed  

Let us conclude this article by giving some more details 
on the canonical map and its image of our surface $S$. 
In what follows, we denote by $W$ the canonical 
image $\varPhi_{|K|} (S)$.  
Moreover, we denote by 
$p_3: \tilde{S}_3 \to S$ the shortest composite of 
quadric transformations such that the variable part 
of $p_3^* |K|$ is free from base points, 
and by $\varphi : \tilde{S}_3 \to W$, the unique morphism 
such that $\varPhi_{|K|} = \varphi \circ p_3^{-1}$. 

First we study the case $\deg \varPhi_{|K|} = 1$. 
In this case, the canonical image $W \subset \mathbb{P}^3$ 
is a sextic surface.  
Recall that for a singularity $(W, x)$ of our surface $W$, 
the fundamental genus of $(W, x)$ is the arithmetic genus 
of its fundamental cycle.
Moreover, since $\varphi : \tilde{S}_3 \to W$ gives the 
minimal desingularization of the canonical image $W$, 
the geometric genus of $(W, x)$ is the dimension of 
the vector space $(R^1 \varphi_* \mathcal{O}_{\tilde{S}_3})_x$, where 
$R^1 \varphi_* \mathcal{O}_{\tilde{S}_3}$ is the first higher direct 
image of the structure sheaf $\mathcal{O}_{\tilde{S}_3}$.  
The following proposition is a comment given to the author 
by Kazuhiro Konno: 

\begin{proposition}   \label{prop:normalityim}
Let $S$ be a minimal algebraic surface as in 
Theorem \ref{thm:structurethm}. Suppose that 
$\deg \varPhi_{|K|} = 1$ and that the canonical system $|K|$ has 
no fixed component. 
Then the canonical image $W = \varPhi_{|K|} (S) \subset \mathbb{P}^3$
is normal. Moreover, if the surface $S$ is sufficiently general, 
then the singularity $(W, x)$ of $W$ is 
a double point with fundamental genus $3$ and geometric genus $6$, 
where $x \in W$ is a point given by 
$(\xi_0 : \eta_0 : \eta_1 : \zeta_0) = (0 : 0: 0: 1)$. 
\end{proposition} 
  
Proof. 
Assume that $|K|$ has no fixed component, as is indeed the case 
for our general $S$ by Proposition \ref{prop:normalform}.  
Then $p_3 : \tilde{S}_3 \to S$ is a blowing up at three 
simple base points of $|K|$. 
Thus for any hyperplane $H \subset \mathbb{P}^3$, 
the arithmetic genus of $W \cap H$ equals that of the 
pullback $\varphi^* (W \cap H) \in |p_3^* (K) -\varepsilon |$, 
where $\varepsilon$ is the sum of the total transforms  
of the three $(-1)$-curves appearing by $p_3 : \tilde{S}_3 \to S$. 
Since the variable part $|p_3^* (K) -\varepsilon |$ of $p_3^* |K|$ 
is free from base points, this together with Bertini's Theorem 
implies that $W$ has at most isolated singularities.  
This however implies that $W$ is normal, 
since the canonical image $W \subset \mathbb{P}^3$ is a hypersurface.   
Note that the local equation at $x$ of $W$ in $\mathbb{P}^3$ is 
analytically in the form $w^2 - f_8 (u, v) = 0$. Thus 
the invariants of the double point $(W, x)$ can be computed 
by the canonical resolution. \qed 

\begin{remark}
From Proposition \ref{prop:normalform}, we see easily that 
the point $x$ is the only singularity of $W$ 
for a sufficiently general $S$. Thus one can compute the 
geometric genus of $(W, x)$ also by writing down the 
Leray spectral sequence of $\varphi : \tilde{S}_3 \to W$ 
and comparing the invariants of 
$W$ and those of $\tilde{S}_3$.   
\end{remark}

Next, we study the 
case $\deg \varPhi_{|K|} = 2$. 
In this case, the canonical image $W \subset \mathbb{P}^3$ is 
a cubic surface.  
We shall describe the branch divisor of the canonical map $\varPhi_{|K|}$.   
For simplicity, we shall do this only for the case where $S$ satisfies 
the following three generality conditions: 
\smallskip

i) the canonical image $W = \varPhi_{|K|} (S)$ is smooth;

ii) the unique member $L \in |L|$ is irreducible; 

iii) the base locus of $|K|$ consists of three distinct points. 

\begin{proposition}   \label{prop:branchdivofphiK}
Let $S$ be a minimal algebraic surface 
as in Theorem \ref{thm:structurethm},    
and $\varphi : \tilde{S}_3 \to W = \varPhi_{|K|} (S)$, 
the morphism such that 
$\varPhi_{|K|} = \varphi \circ p_3^{-1}$.  
Suppose that $\deg \varPhi_{|K|} = 2$ and that $S$ satisfies 
the three conditions above.   
Then the branch divisor $B$ of $\varphi$ 
splits as 
$B = \sum_{i=1}^3 \varGamma_i + B^{\prime}$, where 
$\varGamma_i$'s are 
three coplanar lines in $\mathbb{P}^3$ meeting at one point $x \in W$, 
and $B^{\prime}$ is a member of $|-5 K_W|$ that has an ordinary 
$5$-tuple point at $x$ and such that all other singularities if any are 
negligible ones.     
\end{proposition} 

Proof. 
By the generality conditions, the three base points of the 
canonical system $|K|$ are non-singular points of the unique member $L$. 
Thus if we denote by $\varepsilon$ the divisor such that  
$|K_{\tilde{S}_3}| = p_3^* |K| + \varepsilon$, 
then we have $p_3^* (L) = {p_3}_*^{-1} (L) + \varepsilon$. 
Moreover, the divisor $\varepsilon$ is a sum of three $(-1)$-curves.   
Thus from this together with 
$\varphi^* (-K_W) \sim p_3^* (3L) - \varepsilon$,   
we see that $\varphi_* \varepsilon = \sum_{i=1}^3 \varGamma_i$, 
where $\varGamma_1$, $\varGamma_2$, and $\varGamma_3$ are the  
three lines in $\mathbb{P}^3$ corresponding to the irreducible  
components of the divisor $\varepsilon$. 

Let $R$ and $B = \varphi_* (R)$ be the ramification divisor 
and the branch divisor of 
$\varphi : \tilde{S}_3 \to W$, respectively.  
Then by 
\begin{equation} \label{eql:ramdivphi}
R \sim p_3^* (3L) + \varepsilon - \varphi^*(K_W) 
  \sim 2 \varphi^*(-K_W) + 2 \varepsilon , 
\end{equation}
we have 
$BD 
= (-4K_W + 2 \sum_{i=1}^3 \varGamma_i ) D$ 
for any divisor $D$ on $W$,  
which implies $B \in |-4K_W + 2 \sum_{i=1}^3 \varGamma_i |$. 
Now let us denote by $\tilde{L}_3$ the strict transform 
by ${p_3}$ of the divisor $L$.  
By $\tilde{L}_3 \varphi^* (-K_W) = 0$, we see that 
$\varphi$ contracts $\tilde{L}_3$ to a single point 
$x \in W$, where we have 
$x \in \bigcap_{i=i}^3 \varGamma_i$. 
Moreover, since 
$p_g(S) = h^0 (\mathcal{O}_ W(-K_W )) = 4$ and hence 
$3 \tilde{L}_3 + 2 \varepsilon \in 
|p_3^* (3L) - \varepsilon| = \varphi^*| - K_W|$, 
we obtain a member $\varGamma \in |-K_W|$ such that 
$\varphi^* (\varGamma) = 3\tilde{L}_3 + 2 \varepsilon$ holds. 
Since we have $\varGamma = \sum_{i=1}^3 \varGamma_i$ for this 
$\varGamma$, we see that the three lines   
$\varGamma_1$, $\varGamma_2$, and $\varGamma_3$ 
are coplanar, and that $\varepsilon + \tilde{L}_3 \preceq R$, 
since we have $\varphi (\tilde{L}_3) = \{ x \}$.  
We therefore can put $R = \varepsilon + \tilde{L}_3 + R^{\prime}$, 
where $R^{\prime}$ is a non-negative divisor on $\tilde{S}_3$.  
We put 
$B^{\prime} = \varphi_* (R^{\prime}) \in |- 5K_W|$.       

Now let $\hat{q} : \hat{W} \to W$ be the blowing up at $x$, 
and $\varDelta$, its exceptional divisor. 
Then by $\varphi^* (\varGamma) = 3\tilde{L}_3 + 2 \varepsilon$, 
we obtain $\varphi^* (\varGamma_i) = 
2 \bar{\varGamma}_i + \tilde{L}_3$ for each integer 
$1 \leq i \leq 3$, where $\bar{\varGamma}_i$'s are 
three $(-1)$-curves appearing by $p_3$. 
This implies the liftability of $\varphi : \tilde{S}_3 \to W$ 
to a morphism $\hat{\varphi} : \tilde{S}_3 \to \hat{W}$.
Moreover, we obtain     
$\hat{\varphi}^* (\varDelta) = \tilde{L}_3$.   
Thus $\varDelta$ is not a component of the branch divisor of 
$\hat{\varphi}$, from which we infer  
$\hat{\varphi}_* (R^{\prime}) = \hat{q}_*^{-1} (B^{\prime})$. 
Since we have $\hat{\varphi}_*(R^{\prime}) \varDelta =
R^{\prime} \varepsilon = 5$ by (\ref{eql:ramdivphi}), we see 
from this $\mathrm{ord}_x \, B^{\prime} = 5$. 
But the standard double cover argument implies that 
$\sum_{i=1}^3 \hat{q}_*^{-1} (\varGamma_i) + \hat{q}_*^{-1} (B^{\prime})$ 
has at most negligible singularities. 
Thus the point $x$ is an ordinary $5$-tuple point of $B^{\prime}$, and 
all other singularities of $B^{\prime}$ are negligible ones. 
Finally, the equality $(\sum_{i=1}^3 \varGamma_i \cdot B^{\prime})_x = 15$ 
follows from $\hat{\varphi}^*(\varDelta) = \tilde{L}_3$, since this 
latter implies the absence of singularities lying on $\varDelta$ 
of the divisor  
$\sum_{i=1}^3 \hat{q}_*^{-1} (\varGamma_i) + \hat{q}_*^{-1} (B^{\prime})$. 
\qed 

\begin{remark}
Conversely, a non-singular cubic surface $W \subset \mathbb{P}^3$ and 
a member $B \in |-6K_W|$ having the same properties as in 
Proposition \ref{prop:branchdivofphiK} yields a minimal algebraic 
surface $S$ as in Theorem \ref{thm:structurethm} 
with $\deg \varPhi_{|K|} = 2$. 
Naturally, one easily finds the divisor $L$, guided by the proof 
of the proposition above. 
\end{remark}


\begin{flushright}
\begin{minipage}{25em}
Masaaki Murakami \\
Department of Mathematics and Computer Science\\
Kagoshima University\\
Korimoto 1--21--35 , Kagoshima 890--0065, Japan \\
\texttt{murakami@sci.kagoshima-u.ac.jp} 
\end{minipage}
\end{flushright}

\end{document}